\documentclass[11pt,reqno]{amsart} 
\usepackage[left=2.5cm,right=2.5cm,bottom=3cm]{geometry}
\usepackage{amsfonts}
\usepackage{amsmath}
\usepackage{amssymb}
\usepackage{mathtools}
\usepackage{amsthm}
\usepackage{pgfplots}
\pgfplotsset{compat=newest}
\usepackage[shortlabels]{enumitem}
\usepackage{tikz}   
\usepackage{tikz-cd}  
\usepackage{changepage}
\usetikzlibrary{arrows}
\usepackage{bm}
\usepackage{graphicx}
\usepackage{cases}
\usepackage{appendix}
\usepackage{esint}
\usepackage{hyperref}
\usepackage{mathrsfs}

\numberwithin{equation}{section}

\usepackage{amsmath,amssymb}
\usetikzlibrary{arrows.meta,decorations.markings,calc}
\definecolor{loopone}{HTML}{30876B}
\definecolor{looptwo}{HTML}{BA4A46}
\definecolor{loopthree}{HTML}{2677B8}

\theoremstyle{plain}
\newtheorem{thm}{Theorem}

\newtheorem{lem}[thm]{Lemma}
\newtheorem{prop}[thm]{Proposition}

\numberwithin{thm}{section}

\theoremstyle{definition}

\newtheorem*{question}{Question}
\newtheorem{remark}[thm]{Remark}

\newcommand{\R}{\mathbb{R}}  
\newcommand{\Z}{\mathbb{Z}}

\newcommand{\CC}{\mathbb{C}}   
\newcommand{\Sp}{\mathbb{S}}
  
\newcommand{\Disk}{\mathbb{D}}
\newcommand{\p}{\partial}  
\newcommand{\A}{\mathrm{Area}}

\newcommand{\dif}{\mathrm{d}}

\renewcommand{\Re}{\operatorname{Re}}
\renewcommand{\Im}{\operatorname{Im}}
\newcommand{\Id}{\mathrm{Id}}

\newcommand{\slot}{\mathbin{\mbox{\Large$\lrcorner$}}}
\newcommand{\eps}{\varepsilon}

\usepackage{xcolor}
\usepackage{listings}

\allowdisplaybreaks

\makeatletter
\renewcommand\th@plain{\slshape}
\makeatother

\title
{All stable Schoen--Wolfson cones are minimizing}

\begin{document}

\begin{abstract}
    In the late nineties, R.~Schoen and J.~Wolfson studied the problem of minimizing the area among Lagrangian surfaces in a Riemannian symplectic $4$-manifold. In their parametric approach, they showed that minimizers are smoothly branched-immersed away from a discrete set of points, where there exists a nonflat tangent cone minimizing the area among Lagrangian disk competitors. They also classified the stable singular cones into a family $C_{p,p+1}$ indexed by a positive integer $p$, but left open the question of which of them are globally minimizing. We settle this question by showing that for each $p \geq 1$, $C_{p,p+1} \cap B_1$ minimizes the area among Lagrangian disks spanning its boundary with degree one.
\end{abstract}
\author{Filippo Gaia}
\address{Department of Mathematics, Stanford University,
Stanford, CA 94305, USA}
\email{fgaia@stanford.edu}

\author{Xuanyu Li}
\address{Department of Mathematics, Cornell University,
Ithaca, NY 14853, USA}
\email{xl896@cornell.edu}

\author{Gerard Orriols}
\address{Department of Mathematics, University of Cambridge,
Cambridge CB3 0WA, United Kingdom}
\email{go262@cam.ac.uk}

\date{\today}
\maketitle

\section{Introduction}\label{s: intro}
Lagrangian submanifolds constitute one of the central objects of study in symplectic geometry. If $(M^{2n}, \omega)$ is a symplectic manifold, recall that an immersion $u : \Sigma \to M$ is called \emph{Lagrangian} if $\Sigma$ is $n$-dimensional and $u^*\omega = 0$ on $\Sigma$. This condition is notably flexible, in the sense that a given Lagrangian submanifold admits infinite-dimensional families of perturbations, most notably those induced by Hamiltonian vector fields. In the presence of a compatible Riemannian metric $g$ on $M$, the study of critical points of the $n$-dimensional volume (computed with respect to $g$) among Lagrangian submanifolds becomes particularly rich thanks to this abundance of Lagrangian variations, and sits in the intersection between symplectic and Riemannian geometry.

Historically, the first notion of critical point of the area among Lagrangian submanifolds to be introduced was that of Hamiltonian stationarity (and stability), by Oh \cite{Oh1, Oh2}, and of area minimality in Hamiltonian isotopy classes, always in a smooth setting. A Lagrangian submanifold $\Sigma \subset M$ is \emph{Hamiltonian stationary} if its first variation of the area vanishes for Hamiltonian vector fields $X$, that is, those of the form $X \slot \omega = \dif h$ for a smooth function $h : M \to \R$. Shortly after, Schoen and Wolfson considered the minimization problem in Lagrangian homology classes \cite{SchoenWolfson-lapietra}, and developed an existence and regularity theory of Lagrangian-minimizing surfaces from the classical mapping point of view in their seminal work \cite{Schoen-Wolfson-mapping-problem}.

Their main result is that, in any closed symplectic $4$-manifold $(M, \omega)$ with a compatible metric $g$, every class in $H_2(M, \mathbb{Z})$ annihilated by $[\omega]$ can be represented as the sum of classes $(u_i)_* [\Sigma_i]$, where each $\Sigma_i$ is a closed oriented surface and $u_i : \Sigma_i \to M$ is a weakly conformal Lagrangian Lipschitz branched immersion, which minimizes the Dirichlet energy (and thus the area) among Lagrangian maps homotopic to it. We will refer to such maps as \emph{Lagrangian-minimizing}. In addition, they developed an optimal regularity theory for Lagrangian-minimizing maps, showing smoothness away from a finite set of singularities, which can be classified as either \emph{branch points} or \emph{conical points}; this dichotomy arises from the existence (and uniqueness) of parametrized tangent cones at every point, which is a consequence of the intricate monotonicity formula for Hamiltonian-stationary exact Lagrangian surfaces derived in \cite{Schoen-Wolfson-mapping-problem}.

If $C$ is such a tangent cone, then $C \cap B_1$ is topologically a disk, Lagrangian with respect to the standard symplectic structure $\omega = \dif x_1 \wedge \dif y_1 + \dif x_2 \wedge \dif y_2$ on $\CC^2$, and \emph{Lagrangian-minimizing among disks}, in the sense that no Lagrangian map of a disk whose boundary parametrizes $\partial (C \cap B_1)$ homeomorphically has less area than $C \cap B_1$. In particular, the possible tangent cones are all Hamiltonian-stationary and stable. Note that a cone obtained from a Lagrangian-minimizing map by a parametrized tangent cone procedure is not automatically minimizing among competitors of positive genus, since one cannot use a higher genus surface as the domain of a competitor for the original map.

The classification of Hamiltonian-stationary Lagrangian cones in $\CC^2$, or equivalently of Hamiltonian-stationary Legendrian knots in $\Sp^3$, is elementary (see \cite[Section~7]{Schoen-Wolfson-mapping-problem}): each cone, up to a rotation by an element of $\mathrm{U}(2)$, is a member of the family $C_{p,q}$ indexed by two relatively prime positive integers $p$ and $q$, whose links $\Gamma_{p,q}$ have the explicit forms
\begin{equation}
    \label{eq: link parametrization}
    \Gamma_{p,q} = C_{p,q} \cap \Sp^3 = \gamma_{p,q}(\Sp^1) \subset \CC^2,
    \qquad
    \gamma_{p,q}(e^{i \theta}) = \frac{1}{\sqrt{p+q}} \left( \sqrt{q} e^{i p \theta}, i \sqrt{p} e^{-i q \theta} \right).
\end{equation}
Moreover, $C_{p,q}$ admits the $C^1$ conformal parametrization
\begin{equation}
    \label{eq: cone parametrization}
    u_{p,q} : r e^{i \theta} \in \CC \longmapsto \frac{r^{\sqrt{pq}}}{\sqrt{p+q}} \left( \sqrt{q} e^{i p \theta}, i \sqrt{p} e^{-i q \theta} \right) \in \CC^2.
\end{equation}
The case $p = q = 1$ describes a great circle and corresponds to a flat cone. For noncoprime $p$ and $q$, $u_{p,q}$ is a multiple cover of another cone with relatively prime parameters. The area of $C_{p,q} \cap B_1$ is $\pi \sqrt{pq}$, and for $p,q > 1$, the link $\Gamma_{p,q}$ is a nontrivial torus knot. Since, by a change of orientation, $C_{p,q}$ is $\mathrm{U}(2)$-equivalent to $C_{q,p}$, we may assume that $p \leq q$ for the rest of the article.

A nontrivial calculation, carried out also in \cite[Section~7]{Schoen-Wolfson-mapping-problem}, characterizes the stable cones away from the origin as those with $|p-q| \leq 1$, and the stability can be extended to variations that move the vertex.
Moreover, as a result of Qiu's work~\cite{Qiu-nonorientable}, no non-flat cone can be minimizing among non-orientable surfaces \cite[Proposition~7.8]{Schoen-Wolfson-mapping-problem}.
However, since the work of Schoen and Wolfson from the late nineties, the following issue has remained unresolved:

\begin{question}
Which of the stable Schoen--Wolfson cones $C_{p,p+1}$ do appear as tangent cones to Lagrangian-minimizing maps?
\end{question}

By the stability criterion, no multiple cover of a non-flat cone can be minimizing. This has two immediate consequences: on one hand, in the regularity theory, every non-flat tangent cone to a Lagrangian-minimizing map has multiplicity one (establishing the dichotomy from above). On the other hand, this problem has the peculiar feature that any non-flat minimizing cone would not be so anymore when taken with higher multiplicity; we do not know of any area-minimizing conical current of any dimension and codimension with this property. This fact frustrates any attempt to adapt to the Lagrangian setting any ``geometric'' method to show that a cone is minimizing (e.g.~calibrations, Lawlor projections, etc.).

\subsection{Main result}
The goal of this paper is to show that in fact all such cones are minimizing among Lagrangian disks. Our main result can be stated most naturally with the sharp regularity assumptions, but for that we need to recall two standard definitions. Denoting by $\Disk$ the unit disk in $\CC = \R^2$, we say that a map $u \in H^1(\Disk, \CC^2)$ is weakly Lagrangian if $u^*\omega = 0$ almost everywhere, and denote its (parametrized) area by
\[
    \A(u) = \int_{\Disk}\vert\p_x u\wedge\p_y u\vert\,\dif x\,\dif y.
\]
Recall also \cite{BrezisMironescu} that for a map $F \in H^{1/2}(\Sp^1, \Sp^1)$ one has a well-defined notion of degree, which is an integer, agrees with the usual winding number for continuous maps, and $\deg F \neq 0$ if $F$ is monotone and surjective.

\begin{thm}\label{thm: main theorem}
    Let $u\in H^1(\Disk, \CC^2)$ be a weakly Lagrangian map and assume that $u|_{\p \Disk} = \gamma_{p,p+1} \circ F$ for a map $F \in H^{1/2}(\Sp^1, \Sp^1)$ of degree $1$. Then
    $$\A(u)\geq\pi\sqrt{p(p+1)}=\A(C_{p,p+1}\cap B_1).$$
\end{thm}

This resolves the classification of blowups for the mapping problem for the area of Lagrangian surfaces. In particular, the conjecture that only the first cone $C_{1,2}$ ought to be Lagrangian-minimizing among disks \cite{LeeWang} is not true, although we have found a much more elementary proof of its minimality.
We remark that our method does not settle whether these cones are minimizing with respect to surfaces of higher genus.

\begin{remark}
    There is a corresponding uniqueness statement for every $p$. For $p = 1$, the rigidity follows from our proof (see Remark~\ref{rmk: rigidity p=1}), and for $p > 1$ it can be obtained by a more indirect argument under mild regularity hypotheses on the boundary (see Remark~\ref{rmk: rigidity p>1}).
\end{remark}

\subsection{Previous work}
The existence of at least one singular Lagrangian cone which is minimizing among disks has been known for some time, but by an indirect method. To every conical singular point $x$ of a Lagrangian branched immersion $u : \Sigma^2 \to M^4$ of a closed surface $\Sigma$, one can associate an invariant, its (oriented) Maslov index $\mu(x)$, which describes the winding number of the Lagrangian angle along the link of its tangent cone; it is equal to $p-q$ when the tangent cone is $C_{p,q}$. A computation \cite[Proposition~2.2]{Schoen-Wolfson-mapping-problem} shows that
\begin{equation}
\label{eq: maslov identity}
    \sum_{x \text{ conical}} \mu(x) = \frac{1}{2} \langle c_1(M), u_*[\Sigma] \rangle,
\end{equation}
and therefore, choosing an ambient K\"ahler manifold---such as $(M^4, \omega) = (\Sp^2 \times \Sp^2, \omega_0 \oplus 2 \omega_0)$---with a homology class where $\omega$ vanishes but $c_1(M)$ does not---such as $(2,-1)$---one sees that the left hand side of \eqref{eq: maslov identity} does not vanish, and therefore there must exist at least one minimizing cone. A more refined construction by Wolfson~\cite{Wolfson-nonsmooth} and Micallef--Wolfson~\cite{MicallefWolfson} shows that, even when $c_1(M) = 0$, such as when $M$ is a K3 surface, singularities cannot be avoided.

Before mentioning more details about the proof, let us give some more context and background about the problem. The initial motivation for Schoen and Wolfson to consider the Lagrangian minimizing problem was as an approach to construct Lagrangian minimal surfaces in K\"ahler--Einstein manifolds, which would be Special Lagrangian if the ambient is Calabi--Yau. Indeed, a beautiful computation shows that a smooth closed Lagrangian submanifold in a K\"ahler--Einstein manifold, which is stationary for the area with respect to Lagrangian variations, must be minimal. The existence of these singular minimizing cones, which as mentioned was already known by indirect means, shows that this approach cannot work in full generality, but a better understanding of the variational properties of the cones may still be useful.

Recently there has been some renewed interest in the variational theory of the area for Lagrangian submanifolds. Starting from a simplified monotonicity formula of Rivi\`ere~\cite{Riviere-monotonicity} and the variational framework developed by the same author in \cite{Riviere-Legendrian}, Pigati and Rivi\`ere~\cite{PigatiRiviere} have extended the regularity theory of Schoen and Wolfson to just Hamiltonian-stationary Lagrangian surfaces, interpreted in the sense of parametrized varifolds. The behavior at the boundary has been studied by the first named author in \cite{Gaia2025FreeBoundary}. In higher dimensions, the monotonicity property breaks down for exact Hamiltonian-stationary Lagrangians, and therefore it is unclear whether singularities should be modeled on cones or can be wilder; nonetheless, one can prove a partial regularity result in the setting of Legendrian minimizers \cite{Orriols}. Finally, in the graphical case, one does recover an almost-monotonicity formula which shows that singularities are modeled on special Lagrangian cones, and the lowest dimension in which they exist has recently been established \cite{BhattacharyaOrriolsSkorobogatova}.

To conclude, we also mention two works constructing Hamiltonian-stationary Lagrangian surfaces in $\CC^2$ which have relevant connections to the present paper. In \cite{HeleinRomon-C2}, all smooth Hamiltonian-stationary tori in $\CC^2$ were classified by formulating the relevant PDE as a completely integrable system; a similar link is exploited in Section~\ref{s: a minimizing problem}, where we take advantage of a zero-curvature representation to study minimizers of an auxiliary functional that plays a key role in our proof---see Remark~\ref{rmk: integrable}. In fact, this functional was introduced in \cite{GaiaOrriolsRiviere} to construct Hamiltonian-stationary surfaces with boundary in $\CC^2$ with multiple $(p, p \pm 1)$-conical singularities. The construction therein is variational and relies on a reduction of the stationarity equation to a linear problem involving only one of the complex coordinates. This idea is also central in the present work---see Remark~\ref{rmk: functional} for a more detailed comparison. In light of this work, it would be interesting to know if the surfaces constructed in \cite{GaiaOrriolsRiviere} are also minimizing, at least near the vertex (cf.~\cite{DimlerGaia, HardtSimon}).

\subsection{Structure of the proof}
Given a map $u\in H^1(\Disk,\CC^2)$, which we can write as $u=(u_1,u_2)$ with real and imaginary parts $u_j=f_j+ig_j$, we have
$$u^*\omega=\dif f_1\wedge \dif g_1+\dif f_2\wedge \dif g_2=(J(u_1)+J(u_2))\dif x\wedge \dif y,$$
where $J(u_j)=\det(D u_j)=\p_xf_j\p_yg_j-\p_yf_j\p_xg_j$ is the Jacobian of $u_j$ and $x,y$ are the Euclidean coordinates on $\Disk$. Hence,
\begin{equation}
    \label{eq: cancellation of Jacobian}
    u=(u_1,u_2)\in H^1(\Disk,\CC^2)\text{ is weakly Lagrangian} \quad \iff \quad J(u_1)+J(u_2)=0\text{ a.e.}
\end{equation}
To prove Theorem \ref{thm: main theorem}, we will work with the Dirichlet energy of $u$ instead of its area, with the convention
\[
    E(u) = E(u_1) + E(u_2),
    \qquad
    E(u_j) = \frac{1}{2}\int_{\Disk} \vert\nabla u_j \vert^2 \,\dif x\,\dif y.
\]
By Morrey's $\epsilon$-conformality lemma, it suffices to show the lower bound
\begin{equation}
    \label{eq: main ineq energy}
    E(u) \geq \pi \sqrt{p(p+1)} = \A(C_{p,p+1}\cap B_1).
\end{equation}
In order to prove \eqref{eq: main ineq energy}, we will fix a function $G : \Disk \to \R$ and work instead with the functional
\begin{equation}
\label{eq: functional E+GJ}
E_G(v) = E(v)+\int_{\Disk}GJ(v)\,\dif x\,\dif y,
\qquad
v \in H^1(\Disk, \CC).
\end{equation}
When $u$ is weakly Lagrangian, $E(u) = E_G(u_1) + E_G(u_2)$ thanks to \eqref{eq: cancellation of Jacobian}, and in fact this is the only place where we will use the Lagrangian condition. It is well known that the Dirichlet energy of $v$ is bounded from below by the square of the $H^{1/2}$-seminorm of $v|_{\p \Disk}$, or in Fourier, by a quadratic form in which the $k$-th mode appears with weight $\pi |k|$. The key idea  is that, if one chooses $G(z) = \log |z|^{-1}$ (the Green's function with a pole at the origin) then the weight of the $k$-th mode of $v|_{\p \Disk}$ in the sharp lower bound for $E_G(v)$ is precisely $\pi \sqrt{k(k+1)}$. This is enough to show that \eqref{eq: main ineq energy} holds for maps $u$ whose boundary trace agrees with $\gamma_{p,p+1}$, which contains only the modes $p$ and $-p-1$.

In order to deal with general boundary parametrizations, in the case $p = 1$ it suffices to vary the position of the pole of $G$ so that a balancing condition is satisfied; controlling the remaining Fourier terms using this cancellation proves the minimality---see Section~\ref{s: p=1}. For $p > 1$, the boundary condition is more nonlinear, but as shown in Section~\ref{s: reduction} we can work in a regime in which $u|_{\p \Disk}$ is close to $\gamma_{p,p+1}$. Then one can check that, in a quadratic expansion of the energy, one needs $p$ ``balancing conditions'' to control the terms with an unfavorable sign, which are enforced by choosing $G$ to be the average of $p$ Green's functions with carefully placed poles. Unfortunately, for such a $G$, some of the explicit Fourier calculations from Section~\ref{s: p=1} fail, and some identities must be replaced by the nontrivial inequalities of Section~\ref{s: a minimizing problem}. Finally, the linearization argument and the proof for $p > 1$ are completed in Section~\ref{s: linearization}.

\subsection{Statement on AI usage}
A publicly available model produced an initial proof of Theorem~\ref{thm: main theorem} for $p=1$, which we then digested and simplified. Subsequently, we extended the result to higher $p$ through a long sequence of interactions with AI, and we streamlined the argument, leading to the present manuscript. This article contains no AI-generated text.

\subsection{Acknowledgements}
F.G.~and G.O.~would like to thank Prof.~Tristan Rivière for introducing them to the field of variational problems for Lagrangian surfaces, and for the many discussions on this topic.
X.L.~would like to thank Prof.~Xin Zhou and Prof.~Daniel Stern for their valuable guidance and continuous support. X.L.~also thanks Zhihan Wang for helpful discussions and his interest in this project.
G.O.~is grateful to Prof.~Yang Li for his interest in this work.
F.G. was supported by the Swiss National Science Foundation (SNSF)
through the Postdoc.Mobility grant P500PT\_230344. X.L.~is partly supported by Simons Dissertation Fellowship in Mathematics.

\section{Boundary Fourier series and the cone \texorpdfstring{$C_{1,2}$}{C1,2}}
\label{s: p=1}

To prove Theorem~\ref{thm: main theorem} for $C_{p,p+1}$ we will fix $p$ points $z_1, \dots, z_p \in \Disk$ (the open unit disk in $\CC$) and set
\begin{align}\label{eq:definition-G}
G=\frac{1}{p}\sum_{k=1}^p G_{z_k}, \quad \text{where} \quad G_{z_k}(z) = -\log\frac{\vert z-z_k\vert}{\vert 1-\bar z_kz\vert}\text{ for }z\in \overline{\Disk} .\end{align}
Thus $G$ is the Green's function with poles at $z_1,\dots,z_p$, normalized as follows:
\begin{equation}
    \label{eq: G is green}
    -\Delta G = \frac{2\pi}{p} \sum_{k=1}^p \delta_{z_k}, \qquad
    G > 0 \text{ in } \Disk, \qquad
    G|_{\p\Disk} = 0
    \qquad \text{and} \qquad
    \p_\nu G < 0 \text{ on } \p\Disk.
\end{equation}
In general, $GJ(v)$ is not integrable for $v\in H^1(\Disk,\CC)$, but we have the following result from the theory of integrability by compensation, which is essentially Wente's inequality with the sharp constant as in \cite[Appendix]{BrezisCoron} or \cite{Topping}. Here the dot $\cdot$ is the scalar product in the domain $\R^2$, and our hermitian forms are linear in the first entry and antilinear in the second.

\begin{lem}
    \label{lem: compensated integrability}
    The map
    \[
        (v, w) \in C^\infty_c(\overline{\Disk} \setminus \{ z_1, \dots, z_p \}, \CC)^2
        \longmapsto
        \mathcal{E}_G(v,w) := \int_{\Disk} \nabla v \cdot \nabla \bar{w} \, \mathrm{vol} + i G \, \dif v \wedge \dif \bar{w} \in \CC
    \]
    extends to a bounded hermitian form $H^1(\Disk, \CC) \times H^1(\Disk, \CC) \to \CC$, and moreover $0 \leq E_G(v) = \frac{1}{2} \mathcal{E}_G(v, v) \leq 2 E(v)$ for each $v \in H^1(\Disk, \CC)$. In addition, the functional $E_G$ is conformally invariant.
\end{lem}

\begin{remark}
    \label{rmk: functional}
    The connection between the nonnegative functional $E_G$ and the Schoen--Wolfson cones had already been noticed in joint work by the first and last authors with T.~Rivi\`ere \cite{GaiaOrriolsRiviere}, but interestingly, there it is used in a completely different way. While here it features as a way to relate the energy of the cone maps $u_{p,p+1}$ to their boundary trace in a sharp way, in \cite{GaiaOrriolsRiviere} it arose \emph{a posteriori} from the study of its Euler--Lagrange equations, which are the PDEs enforcing that the Lagrangian angle of a Lagrangian map is a prescribed $\Sp^1$-harmonic map (in this case, the harmonic conjugate of $G$). Moreover the cancellation \eqref{eq: cancellation of Jacobian} is never used in \cite{GaiaOrriolsRiviere}, whereas the Lagrangian angle interpretation plays no role here.
\end{remark}

\begin{proof}[Proof of Lemma~\ref{lem: compensated integrability}]
    The identities $\frac{1}{2} \mathcal{E}_G(v,v) = E_G(v)$ and $\mathcal{E}_G(w,v) = \overline{\mathcal{E}_G(v,w)}$ are immediate. For $p = 1$ the nonnegativity of $E_G$ is precisely the sharp Wente inequality, and for $p > 1$ it follows by averaging the $p$ terms. Conjugating $v$ changes the sign of the Jacobian, and thus,
    \[
        0 \leq E_G(v)
        = 2E(v) - \left(E(\bar{v}) + \int_{\Disk} G J(\bar{v}) \right)
        = 2E(v) - E_G(\bar v) \leq 2 E(v),
    \]
    which yields the boundedness of $\mathcal{E}_G$ in $H^1$. The conformal invariance of both terms is standard.
\end{proof}

Note that the identity $E(u) = E_G(u_1) + E_G(u_2)$ when $u \in H^1(\Disk, \CC^2)$ satisfies \eqref{eq: cancellation of Jacobian} follows from the interpretation of $\int_\Disk G (J(u_1) + J(u_2))$ as the average at $p$ points of the solution $\varphi \in C_0(\overline{\Disk})$ to the Dirichlet problem $-\Delta \varphi = 2\pi (J(u_1) + J(u_2)) = 0$, which is clearly zero.

For the rest of this section we fix $p=1$ and consider $G_0(z) = \log \tfrac{1}{|z|}$, the Green's function with a single pole at $z_1 = 0$.
As explained in the introduction, not only is $E_{G_0}$ nonnegative, but it also controls the Fourier coefficients of the boundary data in a way that is sharp for the components of all the $u_{p,p+1}$.

\begin{lem}\label{lem: minimizing Fourier expansion of energy}
    Given $v \in H^1(\Disk, \CC)$, write $v(e^{i\theta}) = \sum_{n \in \Z} a_n e^{in\theta}$. Then
    \begin{equation}
        \label{eq: sharp lower bound p=1}
        E_{G_0}(v) \geq \pi\sum_{n=-\infty}^{\infty} \sqrt{n(n+1)} \vert a_{n}\vert^2.
    \end{equation}
\end{lem}
\begin{proof}
By approximation (since the right hand side of \eqref{eq: sharp lower bound p=1} is controlled by $\| v \|_{H^{1/2}(\p\Disk)}$) and Lemma~\ref{lem: compensated integrability}, it is enough to assume that $v$ is smooth up to the boundary and vanishes in a neighborhood of the origin.
Use the oriented conformal coordinates $\theta, t$ on $\Disk$ given by $z = e^{-t + i \theta}$, so that $t = G(z)$, and for each $0 \leq t < \infty$ consider the Fourier expansion
\[
    v(e^{-t+i\theta}) = \sum_{n=-\infty}^{\infty} a_{n}(t) e^{in\theta},
    \qquad
    a_n(0) = a_n,
    \quad
    a_n(t) = 0 \text{ for large } t.
\]
Then, integrating by parts and using the AM-GM inequality,
\begin{align*}
    E_{G_0}(v)
    &= \frac{1}{2} \int_{\Disk} |\nabla v|^2 + i G \, \dif v \wedge \dif \bar{v}
    = \frac{1}{2} \int_{\Disk} |\nabla v|^2 - i v \, \dif G \wedge \dif \bar{v} \\
    &= \frac{1}{2} \sum_{n,m \in \Z} \int_0^\infty \int_{0}^{2\pi}
    (a_n'(t) \overline{a_m'(t)}
    + n m a_n(t) \overline{a_m(t)}
    + m a_n(t) \overline{a_m(t)}
    ) e^{i(n-m)\theta}
    \, \dif \theta \, \dif t \\
    &= \pi \sum_{n \in \Z} \int_0^\infty (|a_n'(t)|^2 + n(n+1) |a_n(t)|^2 ) \, \dif t \\
    &\geq -\pi \sum_{n \in \Z} \int_0^\infty 2 \sqrt{n(n+1)} \Re \left[ a_n'(t) \overline{a_n(t)} \right] \, \dif t \\
    &= -\pi \sum_{n \in \Z} \sqrt{n(n+1)} \int_0^\infty \frac{\dif}{\dif t} |a_n(t)|^2 \, \dif t
    = \pi \sum_{n \in \Z} \sqrt{n(n+1)} |a_n(0)|^2.
    \qedhere
\end{align*}
\end{proof}

With this lower bound for $E_{G_0}(u_j)$ in hand, the rest of the proof can be carried out by direct calculations on the Fourier coefficients.

\begin{proof}[Proof of Theorem \ref{thm: main theorem} for $p=1$]
By Morrey's $\epsilon$-conformality lemma \cite[Theorem 1.2]{Morrey1948} (see \cite{FitziWenger2020} for a modern account), for any $\epsilon>0$, we can find a
homeomorphism $\varphi:\overline\Disk\rightarrow \overline\Disk$, such that $$E(u\circ\varphi)\leq\A(u)+\epsilon.$$
    Thus, as mentioned in the introduction, it is enough to prove the energy lower bound \eqref{eq: main ineq energy} for $u = (u_1, u_2)$.
    We may also assume that the centering condition
    \begin{equation}
        \label{eq: centering}
        \int_{0}^{2\pi} u_1(e^{i\theta}) \, \dif \theta = 0
    \end{equation}
    is satisfied.
    Indeed, let $h : \overline{\Disk} \to \CC$ be the harmonic extension of $u_1 |_{\p \Disk}$ and observe that by the Mean Value Formula, the left hand side of \eqref{eq: centering} is $2\pi h(0)$. Since $\sqrt{\frac{3}{2}} h |_{\p\Disk} : \Sp^1 \to \Sp^1$ has degree one in the VMO sense (see~\cite[Section~12.1]{BrezisMironescu}), $h$ must vanish at some interior point $\zeta \in \Disk$. Let $\Phi$ be a M\"obius transformation sending $0$ to $\zeta$, so that $h \circ \Phi(0) = 0$ and $h \circ \Phi$ is the harmonic extension of $u_1 \circ \Phi$ by conformal invariance.
    Then, replacing $u$ by $u \circ \Phi$ preserves $E(u)$ and the hypotheses for $u$, and the left hand side of \eqref{eq: centering} becomes $2\pi h \circ \Phi(0) = 0$.
    
    By hypothesis, we can write
    \begin{equation}
        u|_{\p\Disk} = \gamma_{1,2}\circ F = \left(\sqrt{\frac{2}{3}}F, \frac{i}{\sqrt{3}}\overline{F}^2\right),
        \qquad \text{where }F\in H^{1/2}(\Sp^1,\Sp^1)\text{ has degree }1.
        \label{eq: boundary condition u}
    \end{equation}
    Write the Fourier series of $F$ and $\bar F^2$ as
    \[
        F(e^{i\theta})=\sum_{n=-\infty}^{\infty} a_{n}e^{in\theta}
        \quad \text{and} \quad \bar F(e^{i\theta})^2=\sum_{n=-\infty}^{\infty} b_{n}e^{in\theta}.
    \]
    Using Parseval's theorem, the identities $\| F \|_{L^2}^2 = 2\pi$ and $\deg F = 1$ can be expressed respectively as
    $$\sum_{n=-\infty}^{\infty}\vert a_{n}\vert^2=1 \quad \text{and} \quad \sum_{n=-\infty}^{\infty}n\vert a_{n}\vert^2=1.$$
    Therefore, 
    \[
        \sum_{n=-\infty}^{\infty}\sqrt{n(n+1)}\vert a_{n}\vert^2-\sqrt2
        =\sum_{n=-\infty}^{\infty}\left(\sqrt{n(n+1)}-\sqrt2+1-n\right)\vert a_{n}\vert^2.
    \]
    The only negative coefficient in this sum is in the term $n=0$, but $a_{0} = 0$ by \eqref{eq: centering}.
    The coefficient for $n=1$ vanishes while for other $n$ we have a positive lower bound
    \begin{align*}
        \sqrt{n(n+1)}-\sqrt 2+1-n\geq\begin{cases}
            \alpha := 2-\sqrt2 > 0, &\text{ if }n\leq-1; \\
            \beta := \sqrt 6-1-\sqrt 2 > 0, &\text{ if }n\geq 2.
        \end{cases}
    \end{align*}
    We conclude that 
    \begin{equation}
        \sum_{n=-\infty}^{\infty}\sqrt{n(n+1)}\vert a_{n}\vert^2-\sqrt2
        \geq
        \alpha \sum_{n=-\infty}^{-1}\vert a_{n}\vert^2 + \beta \sum_{n=2}^{\infty}\vert a_{n}\vert^2.
        \label{eq: lower bound for Fourier expansion, positive mode}
    \end{equation}
    Next we look at the corresponding values for $\bar F^2$: the corresponding identities $\lVert\bar F^2\rVert_{L^2}^2 = 2\pi$ and $\deg (\bar F^2) = -2\deg F = -2$ give
    \[
        \sum_{n=-\infty}^{\infty}\vert b_{n}\vert^2=1
        \quad\text{and}\quad
        \sum_{n=-\infty}^{\infty}n\vert b_{n}\vert^2=-2.
    \]
    Hence, adding twice the first sum plus the second one, we get
    \begin{equation}
        \sum_{n=-\infty}^{\infty}\sqrt{n(n+1)}\vert b_{n}\vert^2-\sqrt2
        =\sum_{n=-\infty}^{\infty}\left(\sqrt{n(n+1)}-\sqrt2+2+n\right)\vert b_{n}\vert^2
        \geq - \left(\sqrt 2-1 \right)\vert b_{-1}\vert^2,\label{eq: lower bound for Fourier expansion, negative mode}
    \end{equation}
    where we have used the coarse bound $\sqrt{n(n+1)}-\sqrt2+2+n\geq0$ for $n\ne-1$. Again since $a_{0}=0$, we can compute $b_{-1}$ from $\{ a_{n} \}$ by
    $$b_{-1}=\sum_{n=-\infty}^{\infty} \overline{a_{n}} \, \overline{a_{1-n}}
    =2\sum_{n=-\infty}^{-1} \overline{a_{n}} \, \overline{a_{1-n}}.$$
    By Cauchy--Schwarz, this gives the estimate
    \[
        \vert b_{-1}\vert^2
        \leq 4 x y,
        \qquad\text{where} \quad
        x := \sum_{n\leq -1} \lvert a_{n} \rvert^2
        \quad \text{and} \quad
        y := \sum_{n\geq 2} \lvert a_{n}\rvert^2.
    \]
Putting together Lemma~\ref{lem: minimizing Fourier expansion of energy}, the boundary condition \eqref{eq: boundary condition u}, and the lower bounds \eqref{eq: lower bound for Fourier expansion, positive mode} and \eqref{eq: lower bound for Fourier expansion, negative mode} we get
\begin{align*}
    E(u) - \pi \sqrt{2}
    &= E_{G_0}(u_1) + E_{G_0}(u_2) - \pi \sqrt{2} \\
    &\geq\frac{2\pi}{3} \left( \sum_{n=-\infty}^{\infty}\sqrt{n(n+1)}\vert a_{n}\vert^2 - \sqrt{2} \right) + \frac{\pi}{3} \left( \sum_{n=-\infty}^{\infty}\sqrt{n(n+1)}\vert b_{n}\vert^2 - \sqrt{2} \right) \\
    &\geq \frac{2\pi}{3}(\alpha x + \beta y)-\frac{4\pi}{3}(\sqrt{2}-1)xy.
\end{align*}
It remains to show that this is nonnegative; since $x+y\leq 1$, rearranging gives
\begin{align*}
    \alpha x + \beta y
    \geq (\alpha x + \beta y)(x+y)
    = (\sqrt{\alpha} x - \sqrt{\beta} y)^2 + (\sqrt{\alpha} + \sqrt{\beta})^2 x y
    \geq(\sqrt{\alpha}+\sqrt{\beta})^2xy,
\end{align*}
and the verification that $(\sqrt{\alpha}+\sqrt{\beta})^2 > 0.90 > 0.83 > 2(\sqrt{2}-1)$ concludes the proof.
\end{proof}

\begin{remark}
    \label{rmk: rigidity p=1}
    By analyzing the equality case in Lemma~\ref{lem: minimizing Fourier expansion of energy} and at the end of the proof, we also get directly a rigidity statement. Indeed, if $u \in H^1(\overline{\Disk}, \CC^2)$ is a weakly conformal map satisfying the hypotheses of Theorem~\ref{thm: main theorem} and $\A(u) = \pi \sqrt{2}$, then it follows that $u \equiv u_{1,2}$ up to a M\"obius transformation. 
\end{remark}

\section{Reduction to a near-standard boundary parametrization}\label{s: reduction}

For the rest of the paper, we focus on the proof of Theorem \ref{thm: main theorem} for $p > 1$. In this section we start by showing that, in order to establish the sharp energy lower bound, it suffices to assume that $u|_{\partial \Disk}$ is $C^1$-close to the standard parametrization $\gamma_{p,p+1}$ of $\Gamma_{p,p+1}$ given in \eqref{eq: link parametrization}. Notice that we do not assume anything about the behavior of $u$ in the interior of $\Disk$, so Theorem~\ref{thm: energy lower bound for near standard parametrization} is not a perturbative result.

\begin{thm}\label{thm: energy lower bound for near standard parametrization}
    There exists $\epsilon_p > 0$ with the following property. Let $F\in C^1(\Sp^1,\Sp^1)$ satisfy $\Vert F-\Id_{\Sp^1}\Vert_{C^1}<\epsilon_p$ and $u\in H^1(\Disk,\CC^2)$ be a weakly Lagrangian map with boundary condition
    \[
        u|_{\p\Disk} = \gamma_{p,p+1} \circ F = \frac{1}{\sqrt{2p+1}} \left( \sqrt{p+1}F^p,i\sqrt p\bar F^{p+1} \right).
    \]
    Then, there holds
    $$E(u)\geq\pi\sqrt{p(p+1)}.$$
\end{thm}

Theorem \ref{thm: energy lower bound for near standard parametrization} will be proved in the next two sections; in the rest of this section we explain how to deduce Theorem \ref{thm: main theorem} from it.

\begin{proof}[Proof of Theorem \ref{thm: main theorem} assuming Theorem \ref{thm: energy lower bound for near standard parametrization}]
Let $u\in H^1(\Disk,\CC^2)$ be weakly Lagrangian and assume that 
$$u|_{\p\Disk} = \gamma_{p,p+1} \circ \tilde F = \frac{1}{\sqrt{2p+1}} \left(\sqrt{p+1}\tilde F^p,i\sqrt p \overline{ \tilde F^{p+1} } \right),$$ for some $\tilde F\in H^{1/2}(\mathbb S^1, \mathbb S^1)$ of degree one, not necessarily close to the identity.

Assume by contradiction that $\operatorname{Area}(u)<\pi\sqrt{p(p+1)}$.
Write $ \tilde F(e^{i\theta})=e^{i(\theta+h(\theta))}$, for some $h\in H^{1/2}(\mathbb S^1,\mathbb R)$ (see \cite[Theorem~5.1]{BrezisMironescu}, applied to $e^{-i\theta}\tilde F(e^{i\theta})$).
Let $H\in H^{1}(\mathbb D\setminus \overline{B_{1/2}},\mathbb R)$ such that
\[
\Delta H=0\text{ in }\mathbb D\setminus \overline{B_{1/2}},\quad H\left(\frac{1}{2},\theta\right)=h\left(\theta\right),\quad H\left(1,\theta\right)=0.
\]

Define $\tilde u: \mathbb D\to\CC^2$ as follows:
\begin{align*}
 \tilde u(re^{i\theta})=\begin{cases}
        u(2re^{i\theta})&\text{if }r\leq \frac{1}{2}\\
        \gamma_{p,p+1}(e^{i(\theta+H(r,\theta))})&\text{if }r\geq \frac{1}{2}.
    \end{cases}
\end{align*}
Note that since on the outer annulus $\tilde u$ takes values in $\Gamma_{p,p+1}$, $\tilde u$ is weakly Lagrangian and $\operatorname{Area}(\tilde u)=\operatorname{Area}(u)<\pi\sqrt{p(p+1)}$. Moreover $\tilde u\in H^1(\mathbb D,\CC^2)$, and $\tilde u |_{\p\Disk} =\gamma_{p,p+1}$. Therefore, replacing $u$ by $\tilde u$, we may assume that $u|_{\p\Disk} =\gamma_{p,p+1}$, and we will do so in the following.

Next, by \cite[Theorem~1.4]{FitziWenger2020}, there exists a smooth orientation-preserving diffeomorphism $\Phi: \overline{\mathbb D}\to \overline{\mathbb D}$, conformal near the boundary $\partial \mathbb D$, such that
\[
E(u\circ\Phi)<\pi\sqrt{p(p+1)}.
\]
Note that the trace of $u\circ\Phi$ might differ from $\gamma_{p,p+1}$. 
Since $\Phi^{-1}$ is conformal near $\partial \mathbb D$, its Beltrami coefficient $\mu$ is compactly supported in $\mathbb D$ and $|\mu|$ attains a maximum smaller than $1$. Moreover it satisfies $\mu(\Phi(z)) \overline{\Phi_z(z)} = -\Phi_{\bar z}(z)$.
Extend it by $0$ outside $\Disk$: set
\begin{align*}
    \mu(z)=\begin{cases}
        \dfrac{\p_{\bar z}\Phi^{-1}(z)}{\p_z\Phi^{-1}(z)},&z\in\Disk;\\
        0,&z\in\CC\setminus\Disk.
    \end{cases}
\end{align*}
By the measurable Riemann mapping theorem \cite{AhlforsBers}, there exists a quasi-conformal homeomorphism $\psi:\CC\rightarrow\CC$ fixing $0,1,\infty$ whose Beltrami coefficient equals $\mu$. Then $\psi\circ\Phi:\Disk\rightarrow\psi(\Disk)$ is conformal since
\[
\partial_{\overline z}(\psi\circ\Phi)=(\psi_w\circ\Phi)\Phi_{\overline z}+(\psi_{\overline w}\circ\Phi)\overline{\Phi_z}=(\psi_w\circ\Phi)[\Phi_{\overline z}+(\mu\circ\Phi)\overline {\Phi_z}]=0.
\]

Conformal invariance of the energy in 2 dimensions gives 
\begin{equation}
    E_{\psi(\Disk)}( u\circ \psi^{-1})=E(u\circ\Phi),
    \label{eq:psi-conf-inv}
\end{equation}
where $E_{\Omega}$ denotes the Dirichlet energy on the domain $\Omega\subset\CC$.
Now for $\rho\in (0,1)$, define the following map from $\Disk$ to $\CC^2$.
\begin{align}\label{eq:rescaled-copy}
     u_\rho(re^{i\theta})=\begin{cases}
        \rho^{\sqrt{p(p+1)}} u(re^{i\theta}/\rho),&\text{if }r\leq\rho\\
        u_{p,p+1}(re^{i\theta}),&\text{if }r>\rho.
    \end{cases}
\end{align}
Here $u_{p,p+1}$ denotes the conformal parametrization of $C_{p,p+1}$ from \eqref{eq: cone parametrization}. At $r=\rho$, both traces equal $u_{p,p+1}(\rho e^{i\theta})$, therefore $ u_\rho$ is a weakly Lagrangian map in $H^1(\Disk,\CC^2)$.
We also set $\psi_\rho(z)=\rho\psi(z/\rho)$.

By \eqref{eq:psi-conf-inv} and the conformal invariance of dilations, on the inner disk $B_{\rho}$ we have
\[
    E_{\psi_\rho(B_{\rho})}( u_\rho\circ\psi_\rho^{-1})
    = E_{\psi(\Disk)} \left( \rho^{\sqrt{p(p+1)}} u \circ \psi^{-1} \right)
    =\rho^{2\sqrt{p(p+1)}}E(u\circ\Phi).
\]
Since $\psi_\rho$ is conformal on $\Disk\setminus \overline B_\rho$,

\[
    E_{\psi_\rho(\Disk \setminus B_\rho)}(u_\rho\circ\psi_\rho^{-1})
    = E_{\Disk \setminus B_\rho}(u_\rho)
    = \A(u_{p,p+1}(\Disk \setminus B_\rho))
    = \left( 1-\rho^{2\sqrt{p(p+1)}} \right)\pi\sqrt{p(p+1)}.
\]

As a result,
\begin{align*}
    E_{\psi_\rho(\Disk)}(u_\rho\circ\psi_\rho^{-1})
    =
    \pi\sqrt{p(p+1)}
    - \rho^{2\sqrt{p(p+1)}} \left( \pi\sqrt{p(p+1)} - E(u\circ\Phi) \right)
    <\pi\sqrt{p(p+1)}.
\end{align*}
Since $\psi_\rho$ is holomorphic
in a neighborhood of $\partial \Disk$, the domain $\psi_\rho(\Disk)$ has smooth boundary.
Thus, by the Riemann mapping theorem, there exists a conformal map $R_\rho:\psi_\rho (\Disk)\to \Disk$, which extends to a smooth diffeomorphism $\overline{\psi_\rho(\Disk)}\to\overline\Disk$, such that $R_\rho(0)=0$ and $R_\rho(\psi_\rho(1))=1$.

Set $w_\rho:=u_\rho\circ \psi_\rho^{-1}\circ R_\rho^{-1}$. Then $w_\rho$ is weakly Lagrangian, is in $H^1(\Disk,\CC^2)$ and
\begin{align}
\label{eq:Ew-rho}
E(w_\rho)= E_{\psi_\rho(\Disk)}( u_\rho\circ\psi_\rho^{-1})<\pi\sqrt{p(p+1)}.
\end{align}
Moreover, $w_\rho|_{\p\Disk}=\gamma_{p,p+1}\circ \psi_\rho^{-1}\circ R_\rho^{-1}$. It remains to show that as $\rho$ converges to zero, $F_\rho:=\psi_\rho^{-1}\circ R_\rho^{-1}\vert_{\partial \Disk}$ converges to the identity in $C^1(\partial \Disk)$. To this end, set $\Lambda_\rho=R_\rho\circ \psi_\rho : \overline{\Disk} \to \overline{\Disk}$. Note that since $R_\rho$ is conformal and $\mu_{\psi_\rho}(z)=\mu_\psi(z/\rho)$, the Beltrami coefficients $\mu_{\Lambda_\rho}$ of $\Lambda_\rho$ satisfy $\lVert \mu_{\Lambda_\rho}\rVert_{L^\infty}=\lVert \mu_{\psi}\rVert_{L^\infty}<1$ and are supported in $B_\rho$.
Moreover, for any $\rho\in (0,1)$ $\Lambda_\rho(0)=0$ and $\Lambda_\rho(1)=1$. Therefore, by the compactness properties of normalized quasiconformal maps (see \cite[Section III.C, Theorem 1]{Ahlfors2006}), for any sequence $\rho_j\to 0$ there is a subsequence along which $\Lambda_{\rho_j}$ converges uniformly on $\overline {\mathbb D}$ to a quasiconformal homeomorphism $\Lambda$ of $\Disk$ onto itself with $\Lambda(0)=0$ and $\Lambda(1)=1$.
Since the Beltrami coefficients  of $\Lambda_{\rho_j}$ are supported in smaller and smaller balls around the origin, the limit map $\Lambda$ is holomorphic in $\Disk\setminus\{0\}$, and since it is bounded, it is actually holomorphic on $\Disk$. As such, $\Lambda$ is a conformal automorphism of $\Disk$, fixing $0$ and $1$. Thus $\Lambda$ must correspond to the identity. As this holds for any subsequence $\rho_j\to0$, we conclude that $\Lambda_\rho$ converges uniformly in $\overline \Disk$ to the identity.

Now define the Schwarz reflection at the boundary
\[
\tilde \Lambda_\rho(z)=\begin{cases}
    \Lambda_\rho(z)& \text{if }\lvert z\rvert\leq 1\\
    \overline {\Lambda_\rho(1/\overline z)^{-1}} &\text{if }\lvert z\rvert> 1.
\end{cases}
\]
For $\rho<\frac{1}{2}$, the maps $\tilde \Lambda_\rho$ are holomorphic in $\{1/2<\lvert z\rvert<2\}$, and converge uniformly to the identity there.
By Cauchy estimates,
they converge in $C^1$ in a neighborhood of $\partial \mathbb D$.
Hence, we also have $F_\rho=\Lambda_\rho^{-1}\vert_{\partial\Disk}\to\Id_{\partial \Disk}$ in $C^1(\partial \Disk)$, as $\rho\to 0$.
For sufficiently small $\rho$, however, Theorem \ref{thm: energy lower bound for near standard parametrization} implies that $E(w_\rho)\geq \pi\sqrt{p(p+1)}$, contradicting \eqref{eq:Ew-rho}.
\end{proof}

\begin{remark}
    \label{rmk: rigidity p>1}
        Since our argument to prove Theorem~\ref{thm: main theorem} for $p > 1$ is by contradiction, our method does not establish the rigidity of minimizers directly. Nonetheless, for a regular boundary parametrization, we can obtain a uniqueness statement by an alternative route based only on the regularity theory: suppose that $u \in H^1(\Disk, \CC^2)$ is a weakly conformal map satisfying the hypotheses of Theorem~\ref{thm: main theorem}, and $\A(u) = \A(C_{p,p+1} \cap B_1)$. Assume
        that $u|_{\p\Disk}$ is a $C^1$ embedding. Let $v : \Disk_2 \setminus \Disk \to \CC^2$ be a $C^1$ parametrization of $C_{p,p+1} \cap (B_2 \setminus B_1)$ such that $v|_{\p \Disk} \equiv u|_{\p\Disk}$. The resulting piecewise defined map $\Disk_2 \to \CC^2$ has bounded distortion and therefore, again by the measurable Riemann mapping theorem, it can be factorized as $w \circ \Phi$, where $\Phi : \overline{\Disk_2} \to \overline{\Disk}$ is some quasiconformal homeomorphism, which is conformal in $\Disk$, and $w \in H^1(\Disk, \CC^2)$ is weakly conformal.
    
        Since, by Theorem~\ref{thm: main theorem}, $w$ minimizes the area among Lagrangian disks relative to its boundary, the interior regularity theory of Schoen and Wolfson (see~\cite[Theorem~4.10]{Schoen-Wolfson-mapping-problem}) implies that $w$ is continuous and a smooth immersion in $U \setminus S$, where $U \subset \Disk$ is some neighborhood of $\Phi(\overline{\Disk_1})$ and $S$ is finite. For each point $z_0 \in U \setminus S$ there exists $\rho > 0$ such that $w(B_\rho(z_0))$ is a Hamiltonian-stationary $C^1$ graph in $\CC^2$ and hence is real-analytic \cite[Theorem~1.3]{ChenWarren}. Thus $w(U \setminus S)$ is a connected real-analytic surface containing a portion of $C_{p,p+1}$, and therefore (by density) $w(U)$ is entirely contained in the cone and so is the image of $u$.

        Finally, let $\psi := u_{p,p+1}^{-1} \circ u : \overline{\Disk} \to \CC$, which is also weakly conformal and continuous. Notice that, by the choice of $\Phi$, $u$ is an immersion on $\Omega := \Disk \setminus \Phi^{-1}(S)$. Thus the set $Z := \Omega \cap u^{-1}(0)$ is discrete in $\Omega$, and $\psi$ is holomorphic in $\Omega \setminus Z$. By singularity removal, $\psi$ extends holomorphically to $\Omega$, and then again across $\Phi^{-1}(S)$ to the whole disk.
        Since $\psi$ maps $\p\Disk$ to $\p\Disk$ with degree one, we deduce it must be a M\"obius transformation.
\end{remark}

\section{Minimizers of a relaxed energy functional}\label{s: a minimizing problem}

To prove Theorem \ref{thm: energy lower bound for near standard parametrization} we will need a weaker version of Lemma~\ref{lem: minimizing Fourier expansion of energy} which also holds when $p > 1$. We will obtain it by studying minimizers of a relaxation of $E_G$.
Fix $\tau \in \R$ and $p$ not necessarily distinct points $z_1,...,z_p \in \mathbb D$, and
let $B$ be the Blaschke product
\begin{equation}
    \label{eq: blaschke}
    B(z) = e^{i\tau} \prod_{k=1}^p \frac{z-z_k}{1-\bar z_k z}, \qquad z \in \overline{\Disk}.
\end{equation}
Clearly $B$ is holomorphic in a neighborhood of the closed disk and maps each $z_k$ to zero.
Thus, the Green's function $G$ defined in \eqref{eq:definition-G} becomes
\[
    G(z) = \frac{1}{p} \sum_{k=1}^p G_{z_k}(z) = -\frac{1}{p} \log |B(z)|.
\]
Now let $\varphi:\Sp^1\rightarrow\Sp^1$ be the diffeomorphism defined by
\begin{equation}
    \label{eq: def varphi}
    \varphi^p(z) = B(z), \qquad z \in \p\Disk = \Sp^1.
\end{equation}
This is uniquely determined up to a global phase $e^{2\pi i/p}$.
Denote by $\theta$ the standard coordinate of $\Sp^1$ and by $\eta$ the coordinate on $\Sp^1$ induced by $\varphi$, namely $e^{i\eta}=\varphi(e^{i\theta})$.
Thanks to the Cauchy--Riemann equations for $B$, the two coordinates are related by
$\dif\eta=-\p_{\nu}G\dif\theta$, where $\nu$ is the outer normal to $\p\Disk$. All the norms and Fourier coefficients of functions on $\p\Disk$ in this section are computed with respect to $\eta$.

The goal of this section is to prove the following lower bound, which is effective when the boundary value of a map has most of its energy in the negative Fourier modes. This is a replacement of \eqref{eq: sharp lower bound p=1}, which does not hold when $G$ has more than one pole.

\begin{prop}\label{prop: lower bound through boundary Fourier series}
    For any $F\in H^{1/2}(\Sp^1,\CC)$ and $v \in H^1(\Disk,\CC)$ with $v|_{\p\Disk}=F\circ\varphi$, there holds
    \[
        E_G(v)
        = E(v)+\int_{\Disk}GJ(v) \, \dif x \, \dif y
        \geq-\pi\sum_{n=-\infty}^{\infty}\operatorname{sgn}(n)\sqrt{n(n+1)}\vert F_n\vert^2,
    \]
    where $\{ F_n \}$ are the Fourier coefficients of $F$ with respect to the coordinate $\eta$.
\end{prop}

The idea to prove Proposition~\ref{prop: lower bound through boundary Fourier series} is to minimize $E_G(v)$ over $v$ for a given boundary datum $F$ and study the PDE satisfied by the minimizer. Unfortunately, $E_G$ is not coercive and to guarantee the existence of a minimizer we follow \cite{GaiaOrriolsRiviere} and work instead with the relaxed functionals $E_{sG}$ for $0 < s < 1$. The claim is that
\begin{equation}
    \label{eq: lower bound relaxed}
    E_{sG}(v)
    = E(v) + s\int_{\Disk}GJ(v) \, \dif x \, \dif y
    \geq-\pi\sum_{n=-\infty}^{\infty}\operatorname{sgn}(n)\sqrt{n(n+s)}\vert F_n\vert^2,
\end{equation}
and then Proposition~\ref{prop: lower bound through boundary Fourier series} follows by letting $s \nearrow 1$.

We start with the variational setup:

\begin{lem}\label{lem: minimization DtN}
    For any $s\in(0,1)$ and $F\in H^{1/2}(\p\Disk,\CC)$ there exists a unique minimizer $w \in H^1(\Disk,\CC)$ of $E_{sG}$ under the Dirichlet boundary condition $w |_{\p\Disk} = F \circ \varphi$, and it satisfies the Euler--Lagrange equation
    \begin{equation}
        \label{eq: euler lagrange}
         \dif ({\star} \dif w - i s G \, \dif w ) = 0
         \qquad \text{weakly in } \Disk.
    \end{equation}
    Moreover, let $\Lambda$ be the map that associates to each smooth $F$ the conformally invariant normal derivative $\Lambda F = (-\p_{\nu}G)^{-1} \p_{\nu}w$ of its corresponding minimizer. Then $\Lambda$ extends uniquely to a bounded linear map $\Lambda : H^{1/2}(\Sp^1,\CC)\rightarrow \dot H^{-1/2}(\Sp^1,\CC)$, the weak Dirichlet-to-Neumann operator, and satisfies
    \begin{equation}
        \label{eq: Dirichlet to Neumann}
        \langle \Lambda F, \phi \rangle = \mathcal{E}_{sG}(w, \tilde \phi)
        \qquad \text{ whenever } \tilde{\phi} \in H^1(\Disk, \CC) \text{ is any extension of } \phi  \circ \varphi \in H^{1/2}(\p \Disk, \CC).
    \end{equation}
\end{lem}

Here we denote by $H^{-1/2}(\Sp^1, \CC)$ the conjugate-dual of $H^{1/2}$, meaning the set of bounded anti-linear functionals $H^{1/2}(\Sp^1, \CC) \to \CC$, and by $\dot{H}^{-1/2}(\Sp^1,\CC)\subset H^{-1/2}(\Sp^1,\CC)$ the subspace of those which annihilate the constants.
\begin{proof}
    Recall that $E_{sG}(v) = (1-s) E(v) + s E_G(v)$ and, by Lemma~\ref{lem: compensated integrability},
    $0 \leq E_G(v) \leq 2E(v)$. Thus, for each $s\in(0,1)$,
    \[
        (1-s)E(v) \leq E_{sG}(v) \leq(1+s) E(v),
    \]
    and therefore the quadratic form $E_{sG}$ is positive-definite and coercive, so it admits a unique minimizer $w$ in the affine space $\{ v \in H^1(\Disk, \CC) : v|_{\p\Disk} = F \circ \varphi\}$. A routine computation verifies the Euler--Lagrange equation~\eqref{eq: euler lagrange}.
    
    We define directly the weak Dirichlet-to-Neumann map via \eqref{eq: Dirichlet to Neumann}.
    Since $w$ solves the weak Euler--Lagrange equations, $\mathcal{E}_{sG}(w, \cdot)$ vanishes on $H^1_0(\Disk)$ and thus $\Lambda$ is well-defined; moreover it is clear that $\langle \Lambda F, 1\rangle = \mathcal{E}_{sG}(w,1) = 0$, so $\Lambda F \in \dot{H}^{-1/2}(\Sp^1)$.
    Whenever $F$ is smooth, by standard elliptic regularity, $w$ is smooth around the boundary and then, using \eqref{eq: euler lagrange},
    \begin{align*}
        \langle \Lambda F, \phi \rangle
        &= \int_{\Disk} \nabla w \cdot \nabla \bar{\phi} + i s G \, \dif w \wedge \dif \bar \phi
        = \int_{\Disk} (-{\star} \dif w + i s G \, \dif w ) \wedge \dif \bar \phi \\
        &= \int_{\Disk} \dif (({\star} \dif w - i s G \, \dif w ) \bar \phi) 
        = \int_{\p\Disk} \p_\nu w \, \overline{\phi(e^{i\eta})} \, \dif \theta
        = \int_{\p\Disk} \frac{\p_\nu w}{|\p_\nu G|} \overline{\phi(e^{i\eta})} \, \dif \eta,
    \end{align*}
    so $\Lambda F$ is represented by the function $\psi = |\p_\nu G|^{-1} \p_\nu w$ integrated against $\dif \eta$.
\end{proof}

The core of the proof lies in the following proposition. As one can see from the computations in the proof of Lemma~\ref{lem: minimizing Fourier expansion of energy}, equality holds in the case $p=1$; it is remarkable that for $p>1$ this survives as an inequality.

\begin{prop}\label{prop: upper bound through boundary value}
    Fix $s\in(0,1)$ and let $F\in C^\infty(\Sp^1,\CC)$ and $\psi = \Lambda F \in C^\infty(\Sp^1,\CC)$.
    Then
    $$\sum_{n=-\infty}^{\infty}\frac{\operatorname{sgn}(n)\vert\psi_n\vert^2}{\sqrt{n(n+s)}}\leq\sum_{n=-\infty}^{\infty}\operatorname{sgn}(n)\sqrt{n(n+s)}\vert F_n\vert^2,$$
    where $F_n,\psi_n$ are the $n$-th Fourier coefficients of $F,\psi$ with respect to the coordinate $\eta$.
\end{prop}

We first show how to deduce Proposition~\ref{prop: lower bound through boundary Fourier series} from this. For that we will need an elementary fact from linear algebra. In a finite-dimensional vector space, this follows easily by assuming that an eigenvalue of $T+J$ is negative. As usual, $A \geq B$ means that $A-B$ is positive semidefinite. 
\begin{lem}\label{lem: positivity of linear operators}
    Let $\mathcal{H}$ be a (real or complex) Hilbert space and $T,J:\mathcal{H}\rightarrow\mathcal{H}$ two bounded linear self-adjoint operators. Assume that $T \geq 0$, $J^2=\Id$, and $TJT\leq J$. Then $T+J\geq0$.
\end{lem}
\begin{proof}
    Set $S:=T+J$. Then
    $$ SJS = TJT+2T+J \leq 2T+2J = 2S.$$
    Observe also that $J + \Id = \tfrac{1}{2} (J + \Id)^2 \geq 0$.
    Conjugating with $S$ gives that $SJS + S^2 \geq 0$, and adding up these two inequalities we obtain
    $$S^2+2S = S^2 + SJS + 2S - SJS \geq 0.$$
    Finally, since $S+2\Id=T+J+2\Id\geq\Id$, we can conjugate by its inverse square root and get
    \[
        S = (S+2\Id)^{-1/2} (S^2 + 2S) (S+2\Id)^{-1/2}\geq0.
        \qedhere
    \]
\end{proof}

\begin{proof}[Proof of Proposition \ref{prop: lower bound through boundary Fourier series} assuming Proposition~\ref{prop: upper bound through boundary value}]
    As mentioned above, it is enough to prove \eqref{eq: lower bound relaxed} for a fixed $0 < s < 1$.
    We may assume that $F_0=0$ by subtracting a constant from $v$. Consider the Hilbert space
    $$\mathcal{H} := \lbrace F\in H^{1/2}(\Sp^1,\CC):F_0=0\rbrace\text{ with norm }\Vert F\Vert^2=\langle F,F\rangle=\sum_{n=-\infty}^{\infty}\sqrt{n(n+s)}\vert F_n\vert^2,$$
    which is equivalent to the standard $H^{1/2}$ norm.
    For each $F\in\mathcal{H}$, let $\psi := \Lambda F \in \dot{H}^{-1/2}(\Sp^1,\CC)$ be as in Proposition \ref{prop: upper bound through boundary value} and define
    $$T:\mathcal{H}\rightarrow \mathcal{H}, \qquad F\mapsto\sum_{n \in \Z \setminus \{0\}}(n(n+s))^{-1/2}\psi_ne^{in\eta},$$
    which is bounded by the boundedness of $\Lambda$ into $\dot{H}^{-1/2}$. Also, \eqref{eq: Dirichlet to Neumann} gives
    \begin{equation}
        \label{eq: TF-F}
        \langle TF,F\rangle_{\mathcal{H}}=\sum_{n=-\infty}^{\infty}\psi_n\bar F_n = \frac{1}{2\pi} \langle \Lambda F, F \rangle_{\dot{H}^{-1/2}, H^{1/2}} = \frac{1}{2\pi} \mathcal{E}_{sG}(w,w) = \frac{1}{\pi} E_{sG}(w) \geq 0.
    \end{equation}
    Hence, $T\geq0$. The other operator to consider is $$J:\mathcal{H}\rightarrow\mathcal{H}, \qquad F\mapsto\sum_{n=-\infty}^{\infty}\operatorname{sgn}(n)F_ne^{in\eta}.$$
    By Proposition \ref{prop: upper bound through boundary value} and approximation, there holds 
    \[
        \langle JTF,TF\rangle = \sum_{n \neq 0} \frac{\operatorname{sgn}(n)\vert\psi_n\vert^2}{\sqrt{n(n+s)}} \leq \sum_{n \neq 0} \operatorname{sgn}(n)\sqrt{n(n+s)}\vert F_n\vert^2=\langle JF,F\rangle
        \quad \implies \quad
        TJT\leq J.
    \]
    Now applying Lemma \ref{lem: positivity of linear operators} we deduce that $T+J\geq 0$, and using \eqref{eq: TF-F} we get
    $$0 \leq \pi \langle (T+J)F,F\rangle_{\mathcal{H}} =E(w)+s\int_{\Disk}GJ(w)\dif x\dif y+\pi\sum_{n=-\infty}^{\infty}\operatorname{sgn}(n)\sqrt{n(n+s)}\vert F_n\vert^2,$$
    as desired. For general $v$, \eqref{eq: lower bound relaxed} is clear as $E_{sG}(v) \geq E_{sG}(w)$.
\end{proof}

In the rest of this section we prove Proposition~\ref{prop: upper bound through boundary value}. For that we need to set up some geometric notions.
Fix $0 < s < 1$ and $F \in C^\infty(\Sp^1, \CC)$ with zero average, let $w$ be the minimizer of $E_{sG}$ given by Lemma~\ref{lem: minimization DtN} and let $\psi = \Lambda F$.
Recall that the Euler--Lagrange equation \eqref{eq: euler lagrange} for $w$ asserts the closedness of the $1$-form ${\star} \dif w - i s G \, \dif w$ away from $\{z_k\}_{k=1}^p$. We start by observing that $\alpha \wedge \beta = {\star} \alpha \wedge {\star} \beta$ for any $1$-forms $\alpha$ and $\beta$. This allows us to rewrite \eqref{eq: euler lagrange} as
\begin{equation}
    \label{eq: euler lagrange connection}
    \dif ({\star} \dif w) - i s \, \dif G \wedge \dif w 
    = \dif ({\star} \dif w) - i s \, {\star} \dif G \wedge {\star} \dif w
    = 0,
\end{equation}
which means that ${\star} \dif w$ is closed with respect to the connection $\dif - i s \, {\star} \dif G$. On the other hand, $\dif w$ is automatically closed for the trivial connection. This suggests looking for a family of $1$-forms and connections interpolating between the two, so that each $1$-form is closed with respect to the corresponding connection. One immediately checks that this imposes two conditions, while so far we have described only one free parameter; the other one will be the coefficient of a pure gauge term along $i \dif G$.

For $-\frac{\pi}{2} < t < \frac{\pi}{2}$, set
\begin{equation}
    \label{eq: alpha D}
    \alpha_t := \cos t \, \dif w + \sin t \, {\star} \dif w
    \qquad \text{and} \qquad
    D_t := \dif - i (\sigma_t \, \dif G + \mu_t \, {\star} \dif G),
\end{equation}
and solve for $D_t \alpha_t = 0$. Using \eqref{eq: euler lagrange connection}, this is equivalent to
\[
    s \sin t \, {\star} \dif G \wedge {\star} \dif w
    - (\sigma_t \, \dif G + \mu_t \, {\star} \dif G) \wedge (\cos t \, \dif w + \sin t \, {\star} \dif w)
    = 0.
\]
By expanding we get the system
\[
    \begin{pmatrix}
        \cos t & \sin t \\
        -\sin t & \cos t
    \end{pmatrix}
    \begin{pmatrix} \sigma_t \\ \mu_t \end{pmatrix}
    =
    \begin{pmatrix} s \sin t \\ 0 \end{pmatrix},
\]
hence $\sigma_t = s \cos t \sin t$ and $\mu_t = s \sin^2 t$.
The $1$-forms $\alpha_t$ also satisfy
\begin{equation}
    \label{eq: quadratics alpha}
    |\alpha_t|^2 = |\dif w|^2
    \qquad \text{and} \qquad
    \overline{\alpha_t} \wedge \alpha_t
    = \dif \bar{w} \wedge \dif w
\end{equation}
and, on the boundary,
\begin{equation}
    \label{eq: boundary alpha}
    \alpha_t |_{\p\Disk}
    = (\cos t \, \p_\eta w + \sin t \, |\p_\nu G|^{-1} \p_\nu w ) \, \dif \eta
    = (\cos t \, \p_\eta F + \sin t \, \psi ) \, \dif \eta.
\end{equation}
The connection $D_t$ is clearly hermitian (for the standard metric) and flat, as $G$ is harmonic away from $\{z_k\}_{k=1}^p$. Therefore, parallel transport with respect to $D_t$ induces a local system over $\CC$ with nontrivial monodromy around each of the poles $z_k$. 

\begin{remark}
    \label{rmk: integrable}
    In fact, $D_t$ and $\alpha_t$ can be assembled into a single $(\mathrm{U}(1) \ltimes \CC)$-connection
    \[
        \mathbf{D}_t = \dif + \begin{pmatrix}
            -i (\sigma_t \, \dif G + \mu_t \, {\star} \dif G) & \alpha_t \\
            0 & 0
        \end{pmatrix}
    \]
    whose flatness (for $t \neq 0$) is equivalent to $D_t \alpha_t = 0$ together with the harmonicity of $G$.
    The fact that we get a $1$-parameter family of flat connections associated to the coupled equations for $w$ and $G$ away from the poles is an instance of the phenomenon of zero-curvature representations for integrable PDEs.
    
    Even more is true: our integrable system is locally a reduction of the one considered by H\'elein and Romon in their study of the related Hamiltonian-stationary equation~\cite{HeleinRomon-C2}, where the structure group is $\mathrm{U}(2) \ltimes \CC^2$. This reduction can be seen more clearly at the level of the curvature-free connections using the complex parameter $\lambda = e^{it} \in \CC^\times$. First one shows that the $\mathrm{SU}(2)$ part of the connection can be gauged away, independently of $\lambda$, so the structure group reduces to $\mathrm{U}(1) \ltimes \CC^2$ acting diagonally (see \cite[Remarks~2 and~7]{HeleinRomon-C2}). Then, projecting onto the first complex coordinate (as in~\cite{GaiaOrriolsRiviere} for the Hamiltonian-stationary problem) further reduces the group to $\mathrm{U}(1) \ltimes \CC$ and recovers the connection $\mathbf{D}_t$ (up to a further gauge).
\end{remark}

In this setting we have the following inequality, which for $p=1$ is always an equality:
\begin{lem}
    \label{lem: stokes alpha}
    Let $\alpha$ be a smooth $1$-form on $\overline{\Disk} \setminus \{ z_k \}_{k=1}^p$
    and let $D = \dif - i (\sigma \, \dif G + \mu \, {\star} \dif G)$ be a connection over $\overline{\Disk} \setminus \{ z_k \}_{k=1}^p$, where $\sigma \in \R$ and $0 < \mu < 1$.
    Suppose that $\alpha$ is $D$-closed and $\alpha |_{\p\Disk} = \phi(e^{i\eta}) \, \dif \eta$, where $\phi$ has Fourier coefficients $\{ \phi_n \}$. Assume also that $|\alpha| \in L^2(\Disk)$. Then
    \[
        \int_{\Disk} \Im [\overline{\alpha} \wedge \alpha] \geq 2\pi \sum_{n=-\infty}^{\infty}\frac{\vert \phi_n \vert^2}{n+\mu}.
    \]
\end{lem}

\begin{proof}[Proof of Proposition~\ref{prop: upper bound through boundary value} assuming Lemma~\ref{lem: stokes alpha}]
    Fix a nonzero $t \in \left( -\frac{\pi}{2}, \frac{\pi}{2} \right)$.
    By the Stokes formula and \eqref{eq: quadratics alpha}, we have
    \[
        \int_{\Disk} \Im [\overline{\alpha_t} \wedge \alpha_t]
        = \int_{\Disk} \Im [\dif \bar{w} \wedge \dif w]
        = \int_{\p\Disk} \Im [\bar w\, \dif w]
        = \int_0^{2\pi} \Im [\bar F\p_{\eta}F] \,\dif\eta
        = 2\pi \sum_{n=-\infty}^{\infty}n\vert F_n\vert^2.
    \]
    The first equation in \eqref{eq: quadratics alpha} gives that $|\alpha_t| \in L^2(\Disk)$, thus by Lemma~\ref{lem: stokes alpha} and \eqref{eq: boundary alpha},
    \[
        2\pi \sum_{n=-\infty}^{\infty}n\vert F_n\vert^2
        \geq 2\pi \sum_{n=-\infty}^{\infty}\frac{\lvert i n F_n \cos t + \psi_n \sin t \rvert^2}{n+\mu_t}.
    \]
    Adding this inequality for $t$ and $-t$ (noting that $\mu_t = \mu_{-t} = s \sin^2 t$) cancels the mixed term on the right by the parallelogram identity:
    \[
        \sum_{n=-\infty}^{\infty}n\vert F_n\vert^2
        \geq \sum_{n=-\infty}^{\infty}\frac{n^2 \cos^2 t \, |F_n|^2 + \sin^2 t \, |\psi_n|^2}{n + s \sin^2 t}.
    \]
    Rearranging, this gives
    \[
        \sum_{n=-\infty}^{\infty} \frac{n(n+s) \sin^2 t}{n + s \sin^2 t} |F_n|^2
        = \sum_{n=-\infty}^{\infty} \left( n -\frac{n^2 \cos^2 t}{n + s \sin^2 t} \right) |F_n|^2
        \geq \sum_{n=-\infty}^{\infty} \frac{\sin^2 t}{n + s \sin^2 t} |\psi_n|^2.
    \]
    Now, dividing by $\sin^2 t$ and integrating between $0$ and $\frac{\pi}{2}$, the result follows from the following formula for $n \neq 0$ (recall that $\psi_0 = 0$):
    \[
        \int_0^{\pi/2} \frac{\dif t}{n + s \sin^2 t} = \frac{\pi\operatorname{sgn}(n)}{\sqrt{n(n+s)}}.
        \qedhere
    \]
\end{proof}

    We are left to establish Lemma~\ref{lem: stokes alpha}.
    Since $\alpha$ is $D$-closed and $D$ is flat, $\alpha$ admits locally a $D$-primitive $f$. This can be computed locally using a parallel trivialization: let $H$ locally solve
    \begin{equation}
        \label{eq: parallel H}
        D \left( H^{-1} \right)
        = 0
        \quad \iff \quad
        \dif \log H = -i (\sigma \, \dif G + \mu \, {\star} \dif G).
    \end{equation}
    Then, writing $f = g H^{-1}$, it solves $Df = \alpha$ if and only if $\dif g = H \alpha$, so $g$ can be obtained by integration. By construction, any lift of $f$, $g$ and $H$ along a curve $\gamma$ only depends on the homotopy class of $\gamma$ relative to the endpoints.
    
    To distinguish between the different types of objects, it will be useful to adopt a more invariant language here: $D$ shall be a connection on a $\mathrm{U}(1)$-bundle $\mathcal{L}$, $f$ a section of $\mathcal{L}$, $H$ a section of $\mathcal{L}^{-1}$, $\alpha \in \Omega^1(\mathcal{L})$ and $g$ a complex-valued function.
    
    \begin{lem}
        \label{lem: lift f gamma}
        Let $\gamma : [0, T] \to \overline{\Disk}$ be a piecewise smooth closed curve which encloses $m$ of the poles $z_k$ with winding number one. If $1 \leq m \leq p$, then there exists a unique periodic section $f_\gamma$ of $\mathcal{L}$ along $\gamma$ solving $D f_\gamma = \alpha |_\gamma$.
    
        Moreover, suppose that $H_\gamma$ is a trivializing section of $\gamma^*\mathcal{L}$ and $g_\gamma : [0,T] \to \CC$ is the corresponding trivialization of $f_\gamma$, meaning $f_\gamma(\tau) = g_\gamma(\tau) H_\gamma(\tau)^{-1}$. Then 
        \begin{equation}
            \label{eq: starting condition}
            g_\gamma(0) = \frac{I_\gamma}{\vartheta - 1},
            \qquad \text{where} \qquad
            I_\gamma := \int_0^T H_\gamma \gamma^*\alpha
        \end{equation}
        and $\vartheta := e^{2\pi i m \mu/p}$ is the monodromy of $H_\gamma$.
    \end{lem}
    \begin{proof}
        Solve the equation for $g$ with initial condition $A$: $ g(\tau) = A + \int_0^\tau H_\gamma \gamma^* \alpha $. By \eqref{eq: parallel H},
        the monodromy of $H_\gamma$ is
        \[
            H_\gamma(T) = \vartheta H_\gamma(0),
            \qquad
            \vartheta
            = \exp \left( -i \int_{\gamma} \sigma \, \dif G + \mu \, {\star} \dif G \right)
            = \exp \left( \frac{2\pi i m \mu}{p} \right),
        \]
        and the closing condition for $f_\gamma$ is $g_\gamma(T) = \vartheta A$, with $g_\gamma(T) = A + I_\gamma$. This can always be solved uniquely because $\vartheta \neq 1$, as $m\mu \leq p\mu < p$.
    \end{proof}

    This allows us to define the following functional for loops $\gamma$ as in Lemma \ref{lem: lift f gamma}:
    \[
        \mathcal{Y}(\gamma) := \Im \int_{\gamma} \overline{f_\gamma} \alpha.
    \]
    Choosing $H_\gamma$ and $g_\gamma$ as in the second part of the lemma, with $|H_\gamma|^2 \equiv 1$, gives the useful representation
    \begin{equation}
        \label{eq: Y of g}
        \mathcal{Y}(\gamma)
        = \Im \int_{0}^{T} \overline{H_\gamma f_\gamma(\tau)} H_\gamma \gamma^* \alpha
        = \Im \int_{0}^{T} \overline{g_\gamma(\tau)} \, \dif g_\gamma(\tau).
    \end{equation}
    The crucial property of $\mathcal{Y}$ is its monotonicity with respect to concatenations of loops:
    \begin{lem}
        \label{lem: concatenation}
        Suppose that $\gamma_j : [0, T_j] \to \overline{\Disk}$, $j = 1,2$, are two piecewise smooth loops with a common base point, enclosing with winding number $1$ two disjoint nonempty subsets of poles, whose respective cardinalities are $m_1$ and $m_2$. Let $\gamma_{12} : [0, T_1 + T_2] \to \overline{\Disk}$ be their concatenation. Then
        \[
            \mathcal{Y}(\gamma_1) + \mathcal{Y}(\gamma_2) \geq \mathcal{Y}(\gamma_{12}).
        \]
        \begin{proof}
            Let $H_1$, $H_2$ and $H_{12}$ be the solutions of \eqref{eq: parallel H} along $\gamma_1$, $\gamma_2$ and $\gamma_{12}$, respectively, starting with a normalized common value $H_1(0) = H_2(0) = H_{12}(0)$. Explicitly,
            \[
                H_{12}(\tau) = \begin{cases}
                    H_1(\tau), & 0 \leq \tau \leq T_1 \\
                    \vartheta_1 H_2(\tau - T_1), & T_1 \leq \tau \leq T_1 + T_2,
                \end{cases}
            \]
            where $\vartheta_j$ denotes the monodromy of $H_j$. Write $g_1$, $g_2$ and $g_{12}$ for the corresponding solutions, with
            \[
                A_1 := g_1(0) = \frac{I_1}{\vartheta_1 - 1},
                \qquad
                A_2 := g_2(0) = \frac{I_2}{\vartheta_2 - 1},
                \qquad
                A_{12} := g_{12}(0) = \frac{I_1 + \vartheta_1 I_2}{\vartheta_1 \vartheta_2 - 1}.
            \]
            Using \eqref{eq: Y of g}, $\dif g_1 = \dif g_{12}$ and $\vartheta_1 \dif g_2(\tau) = \dif g_{12}(T_1 + \tau)$, we compute
            \begin{align*}
                &\int_0^{T_1} \overline{g_1} \, \dif g_1
                + \int_0^{T_2} \overline{g_2} \, \dif g_2
                - \int_0^{T_1+T_2} \overline{g_{12}} \, \dif g_{12}
                = \int_0^{T_1} (\overline{g_1} - \overline{g_{12}} ) \, \dif g_{1}
                + \int_{0}^{T_2} (\overline{g_2}  -  \overline{g_{12}(\tau + T_1)} \vartheta_1) \, \dif g_{2}(\tau) \\
                &\quad = (\overline{A_1 - A_{12}}) \int_0^{T_1} \dif g_{1}
                + (\overline{A_2 - \vartheta_1^{-1} g_{12}(T_1)}) \int_{0}^{T_2}  \, \dif g_{2}
                = (\overline{A_1 - A_{12}}) I_1
                + (\overline{A_2 - \vartheta_1^{-1} (A_{12} + I_1)}) I_2
            \end{align*}
            since both $g_{1} - g_{12}$ and $g_2 - \vartheta_1^{-1} g_{12}$ are constant. Now we substitute
            \[
                A_1 - A_{12}
                = A_1 - \frac{(\vartheta_1 - 1) A_1 + \vartheta_1 (\vartheta_2 - 1) A_2}{\vartheta_1 \vartheta_2 - 1}
                = \frac{\vartheta_1 (\vartheta_2 - 1)}{\vartheta_1 \vartheta_2 - 1} (A_1 - A_2)
            \]
            and
            \[
                A_2 - \vartheta_1^{-1} (A_{12} + I_1)
                = A_2 -  \frac{\vartheta_2 (\vartheta_1 - 1) A_1 + (\vartheta_2 - 1) A_2}{\vartheta_1 \vartheta_2 - 1} 
                = \frac{\vartheta_2 (\vartheta_1 - 1)}{\vartheta_1 \vartheta_2 - 1} (A_2 - A_1).
            \]
            Finally, rearranging and using $(m_1 + m_2) \mu < p$ yields
            \begin{align*}
                \mathcal{Y}(\gamma_1) + \mathcal{Y}(\gamma_2) - \mathcal{Y}(\gamma_{12}) 
                &= \Im \left[ \overline{\frac{\vartheta_1 (\vartheta_2 - 1)}{\vartheta_1 \vartheta_2 - 1} (A_1 - A_2)} (\vartheta_1 - 1) A_1 + \overline{\frac{\vartheta_2 (\vartheta_1 - 1)}{\vartheta_1 \vartheta_2 - 1} (A_2 - A_1)} (\vartheta_2 - 1) A_2 \right] \\
                &= \Im \left[ \frac{(1 - \bar\vartheta_1)(1 - \bar\vartheta_2)}{1 - \bar\vartheta_1 \bar\vartheta_2} (\overline{A_1 - A_2})  (A_1 - A_2) \right] \\
                &= \Im \left[ \frac{(1 - e^{-2\pi i m_1 \mu / p})(1 - e^{-2\pi i m_2 \mu / p})}{1 - e^{-2\pi i (m_1+m_2) \mu / p}} \right] \left| A_1 - A_2 \right|^2 \\
                &= 2 \frac{\sin \left( \frac{\pi m_1 \mu}{p} \right) \sin \left( \frac{\pi m_2 \mu}{p} \right)}{\sin \left( \frac{\pi (m_1+m_2) \mu}{p} \right)} \left| A_1 - A_2 \right|^2 \geq 0.
                \qedhere
            \end{align*}
        \end{proof}
    \end{lem}

    \begin{lem}
        \label{lem: Y-boundary}
        If $\alpha |_{\p\Disk} = \phi(e^{i\eta}) \, \dif \eta$, then
        \[
            \mathcal{Y}(\p\Disk) = 2\pi \sum_{n=-\infty}^\infty \frac{|\phi_n|^2}{n+\mu}.
        \]
    \end{lem}
    \begin{proof}
        Write $f \equiv f_{\p\Disk}$ for the periodic solution on the boundary given by Lemma~\ref{lem: lift f gamma}. This satisfies
        \[
            (\p_\eta f + i \mu f) \, \dif \eta
            = (\dif f - i ( \sigma \, \dif G + \mu \, {\star} \dif G) f) |_{\p\Disk}
            = D f |_{\p\Disk}
            = \alpha |_{\p\Disk}
            = \phi(e^{i\eta}) \, \dif \eta.
        \]
        Thus, noting that $0 < \mu < 1$, the Fourier coefficients of $f$ are $f_n = -i (n+\mu)^{-1} \phi_n$, and hence
        \[
            \mathcal{Y}(\p\Disk)
            = \Im \int_0^{2\pi} \bar{f}(e^{i\eta}) \phi(e^{i\eta}) \, \dif \eta
            = 2\pi \Im \sum_{n=-\infty}^\infty \overline{f_n} \phi_n
            = 2\pi \sum_{n=-\infty}^\infty \frac{|\phi_n|^2}{n+\mu}.
            \qedhere
        \]
    \end{proof}

    \begin{lem}
        \label{lem: Y-poles}
        Fix a pole $z_j$ of $G$ and, for $M > 0$ large, denote by $U_M$ the component of $\{ G > M \}$ containing $z_j$. Then there exists a sequence $M_k \to \infty$ such that
        \[
            \mathcal{Y}(\p U_{M_k}) \to 0
            \quad \text{as } \quad k \to \infty.
        \]
    \end{lem}
    \begin{proof}
        Fix a connected level set $\p U_M$ with an arbitrary parametrization, and let $H, g$ and $f$ be as in Lemma~\ref{lem: lift f gamma} with $|H|^2 \equiv 1$.
        Observe that \eqref{eq: starting condition} gives the estimate
        \[
            |g| \leq C \int_{\p U_M} |\alpha| \qquad \text{on } \p U_M.
        \]
        Using the formula \eqref{eq: Y of g} and Cauchy--Schwarz, we can estimate
        \[
            \left| \mathcal{Y}(\p U_M) \right|
            \leq \int_{\p U_M} |\overline{g}| |H| |\alpha|
            \leq C \left( \int_{\p U_M} |\alpha| \right)^2
            \leq C \int_{\p U_M} \frac{|\alpha|^2}{|\nabla G|} \int_{\p U_M} |\nabla G|.
        \]
        Since $|\nabla G| = -\p_\nu G$ on $\p U_M$, using \eqref{eq: G is green} we see that the second factor is just a constant. Integrating from some large $M_0$ to $\infty$ and using coarea, this is
        \[
            \int_{M_0}^\infty \left| \mathcal{Y}(\p U_M) \right| \, \dif M
            \leq C \int_{M_0}^\infty \int_{\{ G = M \}} \frac{|\alpha|^2}{|\nabla G|} \, \dif M
            \leq C \int_{\Disk} |\alpha|^2 < \infty,
        \]
        which yields the conclusion.
    \end{proof}

    \begin{proof}[Proof of Lemma~\ref{lem: stokes alpha}]
        Let $1 \leq p' \leq p$ be the number of distinct poles, and suppose that for $1 \leq j \leq p'$, all the $z_j$ are distinct and each appears with multiplicity $m_j$.
        Fix a point $z_0 \in \Disk$ and choose simple loops $\gamma_{1}, \dots, \gamma_{p'}$ based at $z_0$ with disjoint interiors, each encircling the corresponding pole $z_j$. Let $\Omega$ be the exterior region, and reorder the poles if necessary so that the inner boundary of $\Omega$ is $\gamma := \gamma_{1} * \dots * \gamma_{p'}$, as in Figure~\ref{fig: stokes}.

        \begin{figure}[htbp]
            \centering
            \begin{tikzpicture}[scale=3.5,font=\large]
  \def\ballradius{.11}
  \coordinate (z0) at (-.63,0);
  \draw[black!80,line width=.8pt] (0,0) circle[radius=1];
  \draw[black!80,line width=.8pt,-{Straight Barb[length=1.8mm,width=1.8mm]}]
    (1,0) arc[start angle=0,end angle=19,radius=1];
  \node at (.61,.62) {$\Omega$};
  \foreach \j/\x/\y/\col in {1/-.03/-.58/loopone,2/.49/.02/looptwo,3/-.03/.58/loopthree}{
    \coordinate (z\j) at (\x,\y);
    \pgfmathsetmacro{\axisangle}{atan2(\y,\x+.63)}
    \pgfmathsetmacro{\tangentangle}{acos(\ballradius/veclen(\x+.63,\y))}
    \pgfmathsetmacro{\startangle}{\axisangle+180+\tangentangle}
    \pgfmathsetmacro{\endangle}{\axisangle+540-\tangentangle}
    \coordinate (a\j) at ($(z\j)+(\startangle:\ballradius)$);
    \coordinate (b\j) at ($(z\j)+(\endangle:\ballradius)$);
    \path[fill=\col!9,draw=\col,line width=.85pt]
      (z0)--(a\j)
      arc[start angle=\startangle,end angle=\endangle,radius=\ballradius]
      --cycle;
    \path[\col,line width=.85pt,postaction={decorate},decoration={markings,
      mark=at position .30 with {\arrow{Straight Barb[length=1.6mm,width=1.7mm]}}
    }] (z0)--(a\j);
    \path[\col,line width=.85pt,postaction={decorate},decoration={markings,
      mark=at position .28 with {\arrow{Straight Barb[length=1.6mm,width=1.7mm]}}
    }] (b\j)--(z0);
    \fill (z\j) circle[radius=.01];
    \node[font=\normalsize] at ($(z\j)!.14cm!(z0)$) {$z_{\j}$};
  }
  \node[loopone,below left] at ($(z0)!.55!(a1)$) {$\gamma_1$};
  \node[looptwo,above=5pt] at ($(z0)!.56!(b2)$) {$\gamma_2$};
  \node[loopthree,above left] at ($(z0)!.52!(b3)$) {$\gamma_3$};
  \fill (z0) circle[radius=.015];
  \node[left=5pt] at (z0) {$z_0$};
\end{tikzpicture}
            \caption{The curves $\gamma_j$ enclosing disjoint regions.}
            \label{fig: stokes}
        \end{figure}
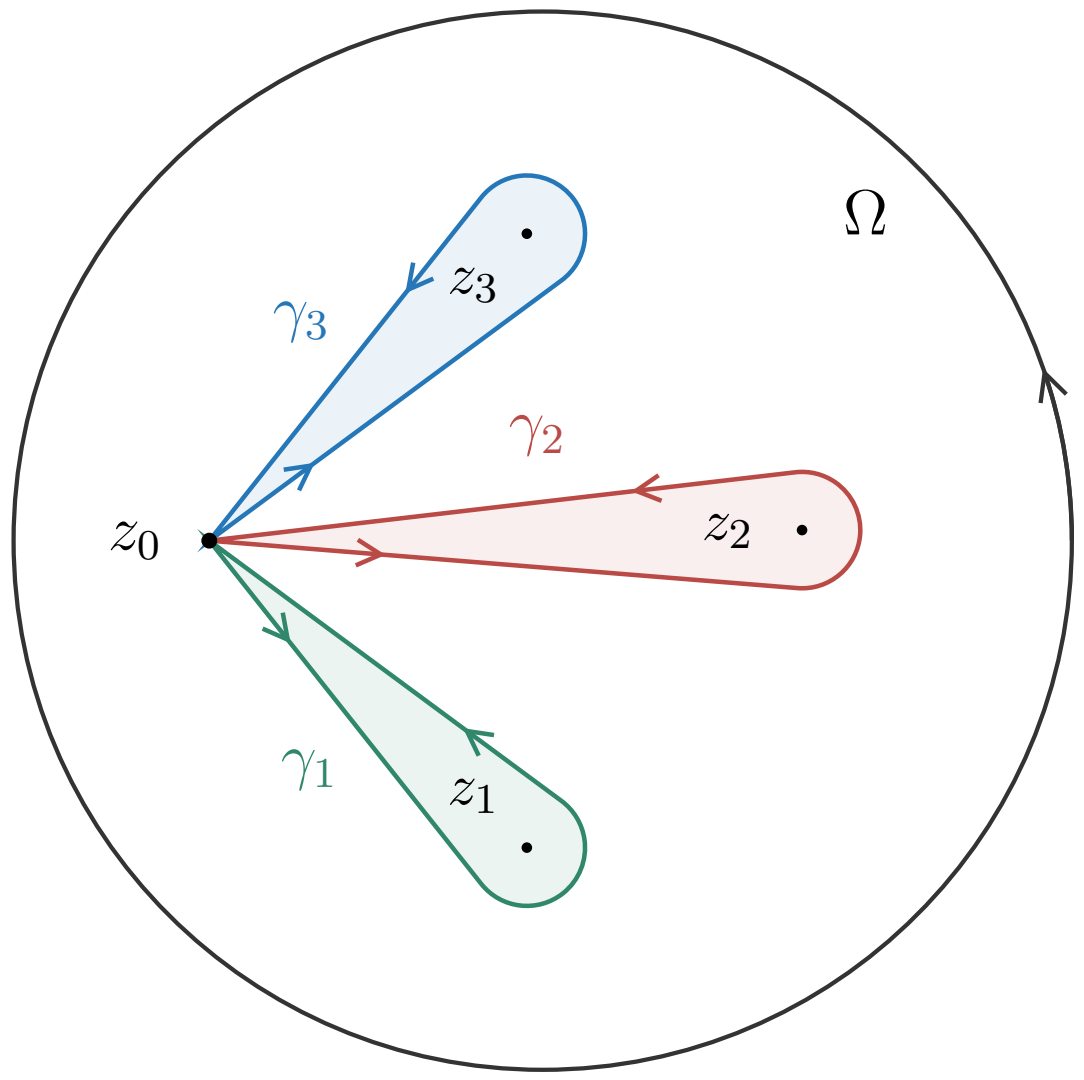
        
        Since $\Omega$ retracts onto $\p \Disk$, by Lemma~\ref{lem: lift f gamma} and homotopy invariance, there exists a unique section $f_0$ of $\mathcal{L}$ over $\Omega$ solving $D f_0 = \alpha$.
        Now a local computation, using a locally defined normalized parallel section $H$ and trivialization $g_0$, shows
        \[
            \overline{\alpha} \wedge \alpha
            = \overline{H \alpha} \wedge H \alpha
            = \overline{\dif g_0} \wedge H \alpha
            = \dif (\overline{g_0} H \alpha)
            = \dif (\overline{Hf_0} H \alpha)
            = \dif (\overline{f_0} \alpha)
        \]
        (one can also argue globally using the unitarity of $D$). Hence the Stokes theorem gives
        \[
            \Im \int_{\Omega} \overline{\alpha} \wedge \alpha
            = \Im \int_{\Omega} \dif (\overline{f_0} \alpha)
            = \Im \int_{\p \Disk} \overline{f_0} \alpha - \Im \int_{\gamma} \overline{f_0} \alpha
            = \mathcal{Y}(\p \Disk) - \mathcal{Y}(\gamma)
        \]
        and thus, applying Lemma~\ref{lem: concatenation} repeatedly,
        \begin{equation}
            \label{eq: stokes exterior}
            \mathcal{Y}(\p \Disk) -\int_{\Omega} \Im[ \overline{\alpha} \wedge \alpha ]
            = \mathcal{Y}(\gamma)
            \leq \mathcal{Y}(\gamma_{1}) + \dots + \mathcal{Y}(\gamma_{p'})
        \end{equation}
        Now for each $j$ choose a sequence $M_k \to \infty$ such that the conclusion of Lemma~\ref{lem: Y-poles} holds, and apply the Stokes theorem again, denoting by $\Omega_j$ the region enclosed by $\gamma_j$:
        \begin{equation}
            \label{eq: stokes interior}
            \mathcal{Y}(\gamma_{j})
            = \int_{\Omega_j \setminus U_{M_k}} \Im [\overline{\alpha} \wedge \alpha] + \mathcal{Y}(\p U_{M_k})
            \xrightarrow{k \to \infty} 
            \int_{\Omega_j} \Im [\overline{\alpha} \wedge \alpha].
        \end{equation}
        This is justified by the integrability hypothesis. Adding \eqref{eq: stokes exterior} and \eqref{eq: stokes interior}, we conclude by Lemma~\ref{lem: Y-boundary}.
    \end{proof}

\section{Proof of the theorem with nearly standard boundary parametrization}
\label{s: linearization}

In the previous section we have shown, in particular, a sharp lower bound on $E_G$ for maps whose boundary data is exactly the second component of $\gamma_{p,p+1} \circ \varphi$. Here we start by showing the analogous lower bound for the first component, which is much easier because it admits an explicit minimizer for every $p$.
In fact we prove a stronger result, which is valid for more general boundary data and retains a gain from the difference on the boundary. Recall the definitions \eqref{eq: blaschke}--\eqref{eq: def varphi} of the Blaschke product $B$ and the boundary diffeomorphism $\varphi$ for fixed parameters $z_1, \dots, z_p \in \Disk$ and $\tau \in \R$.

\begin{prop}
    \label{prop: lower bound positive part}
    The function $w_B := |B|^{\sqrt{(p+1)/p}-1} B$ minimizes $E_G$ for its boundary values $w_B |_{\p\Disk} = \varphi^p$. Moreover, for any $v \in H^{1}(\Disk, \CC)$,
    \begin{align*}
        E_G(v)
        = \sqrt{p(p+1)} \left( \pi - \int_{\p\Disk} \Re [1 - \varphi^{p} \bar{v}] \, \dif\eta \right) + E_G(v - w_B).
    \end{align*}
\end{prop}
\begin{proof}
    Set $\alpha := \sqrt{\frac{p+1}{p}} > 1$ and write $w_B = w_0 \circ B$, where $w_0(z) = |z|^{\alpha-1} z$. It is clear that $w_0 \in C^{1}(\overline{\Disk})$, and it is smooth away from the origin, so $w_B \in C^1(\overline{\Disk}) \cap C^\infty(\Disk \setminus \{ z_k \}_{k=1}^p)$. We first show that $w_B$ solves the Euler--Lagrange equation \eqref{eq: euler lagrange} (with $s=1$) pointwise away from $\{ z_k \}_{k=1}^p$. To do that, we use the fact that pulling back by $B$ commutes with $\star$ on $1$-forms, as $B$ is holomorphic, and $G = p^{-1} G_0 \circ B$, where $G_0(z) = -\log |z|$:
     \[
         \dif ({\star} \dif w_B - i G \, \dif w_B )
         = \dif ({\star} \dif (B^* w_0) - i G \, \dif (B^* w_0) )
         = B^* \dif ({\star} \dif w_0 - i p^{-1} G_0 \,\dif w_0 ),
     \]
     so it suffices to check that $-p \, {\star} \dif w_0 + i G_0 \,\dif w_0$ is closed. In the conformal coordinates $\theta, t$ used in the proof of Lemma~\ref{lem: minimizing Fourier expansion of energy}, $t = G_0$ and $w_0 = e^{-\alpha t + i\theta}$, thus
     \[
         p \, {\star} \dif w_0 - i G_0 \, \dif w_0
         = ((\alpha p + t) \dif \theta + i (p + \alpha t) \dif t) w_0
     \]
     and
     \[
         \dif (p \, {\star} \dif w_0 - i G_0 \, \dif w_0)
         = ( -1 + \alpha (\alpha p + t) - (p + \alpha t) ) w_0 \dif \theta \wedge \dif t
         = ( \alpha^2 p - p - 1 ) \, w_0 \dif \theta \wedge \dif t
         = 0.
     \]

    Next we see that, by a standard cutoff argument, $w_B$ is a weak solution of \eqref{eq: euler lagrange} on the whole $\Disk$: given a test function $\phi \in C^{\infty}_c(\Disk, \CC)$, if $\chi_\eps \in C^\infty_c(\cup_k B_{2\eps}(z_k))$ satisfies $\chi_\eps \equiv 1$ on $\cup_k B_\eps(z_k)$ and $|\dif \chi_\eps| \leq 2 \eps^{-1}$,
    \[
        |\mathcal{E}_G(w_B, \phi)|
        \leq |\mathcal{E}_G(w_B, \chi_\eps \phi)| + |\mathcal{E}_G(w_B, (1 - \chi_\eps) \phi)|
        \leq C  \eps \lvert \log \eps \rvert \| w_B \|_{C^{1}} \| \phi \|_{C^{1}}
        \xrightarrow{\eps \to 0} 0.
    \]
    
    Next we compute $E_G(w_B)$ using \eqref{eq: Dirichlet to Neumann}, which is also valid for $s=1$:
    write $B$ in polar coordinates as $B = e^{-pG + i \Theta}$, so that $w_B = e^{-\alpha p G + i \Theta}$, where $\Theta$ is $\R / 2\pi\Z$-valued. The Cauchy--Riemann equations for $B$ give that $\p_\nu \Theta = 0$ on $\p\Disk$, thus $\p_\nu w_B = (-\p_\nu G) \alpha p B$, and therefore, since $w_B$ is smooth near $\p\Disk$,
    \[
        E_G(w_B)
        = \frac{1}{2} \mathcal{E}_G(w_B, w_B)
        = \frac{1}{2} \langle |\p_\nu G|^{-1} \p_\nu w_B, w_B \rangle_{L^2(\p\Disk, \dif \eta)}
        = \frac{1}{2} \int_{0}^{2\pi} \alpha p B \bar{B} \, \dif \eta
        = \pi \sqrt{p(p+1)}.
    \]
    Now given any $v$ as in the statement, using again \eqref{eq: Dirichlet to Neumann} we obtain
    \begin{align*}
        E_G(v)
        &= E_G(w_B) + \Re \mathcal{E}_G(w_B, v - w_B) + E_G(v - w_B)\\
        &= \pi\sqrt{p(p+1)} + \Re \int_{\p\Disk} \alpha p \varphi^p (\bar v - \bar{\varphi}^p) \, \dif\eta + E_G(v - w_B)\\
        &= \pi\sqrt{p(p+1)} - \sqrt{p(p+1)} \int_{\p\Disk} \Re [1 - \varphi^{p} \bar{v}] \, \dif\eta + E_G(v - w_B).
    \end{align*}
    The minimality of $w_B$ is immediate by the nonnegativity of $E_G$.
\end{proof}

The next technical lemma tells us that choosing the position of the $p$ poles $z_1,...,z_p$ allows us to cancel the first $p$ Fourier coefficients of the phase difference between the given boundary parametrization and $\varphi$. This plays an analogous role to the centering condition \eqref{eq: centering} in the proof of the minimality of $C_{1,2}$.
\begin{lem}\label{lem: normalized boundary parametrization}
    For any $C^1$ diffeomorphism $f:\Sp^1\rightarrow\Sp^1$ such that $\Vert f-\Id_{\Sp^1}\Vert_{C^1}$ is small enough, we can find $z_1,\dots,z_p\in\Disk$ and $\tau\in \mathbb R$ with the following property.
    For the diffeomorphism $\varphi:\Sp^1\rightarrow\Sp^1$ with $\varphi^p(z)=e^{i\tau}\prod_{k=1}^p(z-z_k)/(1 - \overline{z_k} z)$ (choosing the branch close to the identity) and the new coordinate $\eta$ on $\Sp^1$ as in Section \ref{s: a minimizing problem}, there holds
        $$f\circ\varphi^{-1}(e^{i\eta})=e^{i(\eta+g(\eta))},$$
        where $g\in C^1(\R/2\pi\Z)$ and $g_n=0$ for all $\vert n\vert\leq p$. Moreover, for all $\delta>0$, there exists $\epsilon=\epsilon(p,\delta)$ such that if $\Vert f-\Id_{\Sp^1}\Vert_{C^1}<\epsilon$, then $\max_k\vert z_k\vert +\lVert \varphi-\Id_{\Sp^1}\rVert_{C^1}+\Vert g\Vert_{C^1}<\delta$.
\end{lem}
\begin{proof}
    Let us write $f(e^{i\theta})=e^{i(\theta+d(\theta))}$ for a $2\pi$-periodic function $d$. For $z_1,\dots,z_p\in B_{\delta}$ with $\delta$ small, write
    \[
        \prod_{k=1}^p(z-z_k)=z^p+a_{1}z^{p-1}+\cdots+a_p
        \quad \Longrightarrow \quad
        \prod_{k=1}^p(1 - \overline{z_k} z)= 1 + \overline{a_{1}} z + \cdots + \overline{a_p} z^p.
    \]
    We solve for $a_1,\dots a_{p}$ instead of $z_1,\dots, z_p$. One can write $\varphi$ explicitly as $\varphi(e^{i\theta})=e^{i(\theta+L(\theta))}$, corresponding to $\eta(\theta)=\theta+L(\theta)$, where
    $$L(\theta)=\frac{\tau}{p}-\frac{i}{p}\left(\log\left(1+\sum_{k=1}^pa_ke^{-ik\theta}\right)-\log\left(1+\sum_{k=1}^p\overline{a_k}e^{ik\theta}\right)\right)\in\R.$$
    Here we only need to use a branch of $\log$ near $1$ when we have selected $\vert a_k\vert$ small.

    The new function $g(\eta(\theta))=d(\theta)-L(\theta)$, whose Fourier coefficients are
    $$g_n=\frac{1}{2\pi}\int_0^{2\pi}g(\eta)e^{-in\eta}\dif\eta=\frac{1}{2\pi}\int_0^{2\pi}(d(\theta)-L(\theta))e^{-in(\theta+L(\theta))}(1+L'(\theta))\dif\theta.$$
    Since $g$ is a real valued function, they satisfy $g_n=\overline {g_{-n}}$, so we only need to guarantee $g_n=0$ for $0\leq n\leq p$. Letting $\mathscr{G}(d, \tau, a_1, \dots, a_{p})$ be the tuple $(g_n)_{n=0}^p$ of the first $p+1$ Fourier coefficients of $g(\eta) = d(\theta(\eta)) - L(\theta(\eta))$, we have defined a smooth map (note that $g_0 \in \R$)
    $$ \mathscr{G} : C^1(\R/2\pi\Z)\times\R\times\CC^p\rightarrow\R\times\CC^p.$$
    Clearly $\mathscr{G}(0)=0$, and to find its differential at $0$ we first linearize $L(\theta)$
    at $(0,...,0)$,
    $$ D L_0(s,v_1,...,v_p)(\theta) =\frac{s}{p}-\frac{i}{p}\sum_{k=1}^p(v_ke^{-ik\theta}-\bar v_k e^{ik\theta}),$$
    and then compute for $X\in C^1(\R/2\pi),s\in\R$ and $v=(v_1,\dots,v_p)\in\CC^p$:
    \[
        D\mathscr{G}_0(X,s,v)
        = \left(\frac{1}{2\pi}\int_0^{2\pi}(X- D L_0(s,v))e^{-in\theta} \, \dif \theta \right)_{n=0}^p
        = \left( X_0 - \frac{s}{p}, \left( X_n - i \frac{\bar v_n}{p} \right)_{n=1}^p \right),
    \]
    where $\{X_n\}$ are the Fourier coefficients of $X$ with respect to $\theta$. In particular, $D\mathscr{G}_0(0,\cdot,\cdot):\R\times\CC^p\rightarrow\R\times\CC^p$ is obviously invertible. By the implicit function theorem, we can find neighborhoods $U$ of $0$ in $C^1(\R/2\pi\Z)$ and $V$ of $0$ in $\R\times\CC^p$, and a smooth map $\mathscr{A}:U\rightarrow V$ with $\mathscr{A}(0)=0$, such that $\mathscr{G}(d,\mathscr{A}(d))=0$ for all $d\in U$. Since $\max \{ |z_1|, \dots, |z_p| \}$ depends continuously on the coefficients $\mathscr{A}(d)$, and $\lVert \varphi - \Id \rVert_{C^1(\Sp^1)}$  and $\|g\|_{C^1}$ depend continuously on the data,
    the smallness of the norms follows from the continuity of $\mathscr{A}$.
\end{proof}

We are now ready to prove Theorem \ref{thm: energy lower bound for near standard parametrization}.
\begin{proof}[Proof of Theorem \ref{thm: energy lower bound for near standard parametrization}]
    Take $z_1,\dots,z_p,\tau$ and the corresponding $\varphi,g$ given by Lemma \ref{lem: normalized boundary parametrization} (applied with $f=F$), and write
    $$F\circ\varphi^{-1}(e^{i\eta})=e^{i(\eta+g(\eta))},
    \qquad \text{where }g_n=0\text{ for }\vert n\vert\leq p.$$
    Let $\delta>0$ be a small parameter, to be chosen later. 
    Once $\delta$ is fixed, the choice $\epsilon_{p}= \epsilon(p,\delta) > 0$
    will yield
    \begin{align}\label{eq:choice-of-delta}
        \max_k\vert z_k\vert +\lVert \varphi-\Id_{\Sp^1}\rVert_{C^1}+\Vert g\Vert_{C^1}<\delta.
    \end{align}
    Hereafter we use the normalization $u = (2p+1)^{-1/2}(\sqrt{p+1}u_1,i\sqrt p u_2)$.
    Since $u$ is Lagrangian, there holds
    \begin{align*}
        (p+1)J(u_1)+pJ(u_2)=0.
    \end{align*}
    Therefore, for $G$ as in \eqref{eq:definition-G},
    \begin{equation}
        \label{eq: energy u1 u2} 
        E(u)=\frac{1}{2p+1}\left((p+1) E_G(u_1) + p E_G(u_2)\right).
    \end{equation}
    Let
    \[
        F^+(e^{i\eta})=e^{ip\eta}(e^{ipg(\eta)}-1), \qquad F^-(e^{i\eta})=e^{-i(p+1)\eta}(e^{-i(p+1)g(\eta)}-1),
    \]
    so that
    \[
        u_1(\varphi^{-1}(e^{i\eta}))=e^{ip\eta}+F^+(e^{i\eta}), \qquad u_2(\varphi^{-1}(e^{i\eta}))=e^{-i(p+1)\eta}+F^-(e^{i\eta}).
    \]
    
    We will now obtain lower bounds for $E_G(u_1)$ and $E_G(u_2)$. Since $u_1 (\eta)= e^{ip(\eta + g(\eta))}$ on $\p\Disk$, $\Re [1 - \varphi^p \overline{u_1}] = \Re [1 - e^{-ipg} ] = \frac{1}{2}|1 - e^{ipg}|^2 $, thus Proposition~\ref{prop: lower bound positive part} gives
    \[
        E_G(u_1)
        = \pi\sqrt{p(p+1)}\left(1-\frac{1}{2\pi}\Vert e^{ipg}-1\Vert^2_{L^2(\dif\eta)}\right) + E_G(u_1 - w_B).
    \]
    Since $F^+ \circ \varphi$ is the boundary trace of $u_1 - w_B$, using Lemma \ref{lem: lower bound for perturbed single pole problem} below, there holds
    \begin{align*}
        E_G(u_1 - w_B)=
        \frac{1}{p} \sum_{k=1}^p E_{G_{z_k}}(u_1 - w_B)
        \geq\pi\sum_{n=-\infty}^{\infty}\sqrt{n(n+1)}\vert F^+_n\vert^2-C\delta^{1/2}\Vert F^+\Vert_{H^{1/2}(\partial \mathbb D)}^2
    \end{align*}
    provided that $\delta$ is sufficiently small,
    where the Fourier coefficients and the norms are taken in the $\eta$ coordinate.
    Combining these estimates we obtain
    \begin{align}
        E_G(u_1)
        &\geq \pi \sqrt{p(p+1)} \left(1-\frac{1}{2\pi}\Vert e^{ipg}-1\Vert^2_{L^2}\right)+\pi\sum_{n=-\infty}^{\infty}\sqrt{n(n+1)}\vert F^+_n\vert^2
        -C\delta^{1/2} \Vert F^+\Vert_{H^{1/2}}^2.\label{eq: lower bound for positive part}
    \end{align}
    The analogous bound for $u_2$ comes from Proposition \ref{prop: lower bound through boundary Fourier series}, which gives
    \begin{equation}
        \begin{aligned}
            E_G(u_2)
            &\geq-\pi\sum_{n=-\infty}^{\infty}\operatorname{sgn}(n)\sqrt{n(n+1)}\vert(
            u_2\vert_{\partial \mathbb D}\circ\varphi^{-1})_n\vert^2\\
            &=\pi\sqrt{p(p+1)}\left(1-\frac{1}{2\pi}\Vert e^{-i(p+1)g}-1\Vert^2_{L^2}\right)-\pi\sum_{n=-\infty}^{\infty}\operatorname{sgn}(n)\sqrt{n(n+1)}\vert F^-_n\vert^2.
        \end{aligned}\label{eq: lower bound for negative part}
    \end{equation}
    Here we are using that $(u_2\vert_{\partial \mathbb D}\circ\varphi^{-1})_n = \delta_{n,-p-1}+ F^-_n$ for any $n\in \mathbb Z$, and the computation
    \begin{align*}
        \left\lvert 1 + F^-_{-p-1}\right\rvert^2
        &=1+2\Re F^-_{-p-1}+\lvert F^-_{-p-1}\rvert^2\\
        &=1+\frac{1}{\pi}\Re\int_0^{2\pi}(e^{-i(p+1)g(\eta)}-1)d\eta+\lvert F^-_{-p-1}\rvert^2\\
        &=1-\frac{1}{2\pi}\int_0^{2\pi}\vert e^{-i(p+1)g(\eta)}-1\rvert ^2d\eta+\lvert F^-_{-p-1}\rvert^2.
    \end{align*}
    In what follows, we will need two auxiliary estimates:
    \begin{align}
        \label{eq: linearization L2 norm}
         \Vert e^{ing}-1\Vert_{L^2}^2&\leq n^2\Vert g\Vert_{L^2}^2 
        \\
        \label{eq: estimate F+}
        \Vert F^+\Vert_{H^{1/2}}
        &\leq C_p\Vert g\Vert_{H^{1/2}}.
    \end{align}
    The first one follows from the pointwise bound $\lvert e^{it}-1\rvert = 2\,\lvert\sin(t/2)\rvert\leq \lvert t\rvert$,
    while for the second one we use the fact that multiplication by $e^{ip\eta}$ is continuous on $H^{1/2}$, and the estimate $\lvert e^{ipa}-e^{ipb}\rvert\leq p\lvert a - b \rvert$.

    Before combining \eqref{eq: lower bound for positive part} and \eqref{eq: lower bound for negative part}, notice that \eqref{eq: linearization L2 norm} gives
    \begin{align*}
        &\frac{1}{2\pi} \frac{p}{2p+1} \Vert e^{-i(p+1)g}-1\Vert^2_{L^2} + \frac{1}{2\pi} \frac{p+1}{2p+1} \lVert e^{ipg}-1 \rVert^2_{L^2} \\
        &\quad \leq
        \frac{1}{2\pi} \left( \frac{p}{2p+1} (p+1)^2 + \frac{p+1}{2p+1} p^2 \right) \|g\|_{L^2}^2=
        p(p+1) \sum_{n=-\infty}^\infty |g_n|^2.
    \end{align*}
    Hence, inserting \eqref{eq: lower bound for positive part} and \eqref{eq: lower bound for negative part} into \eqref{eq: energy u1 u2} and applying \eqref{eq: estimate F+}, we get
    \begin{equation}
        \label{eq: aux lower bound E}
        \begin{aligned}
            E(u)
            &\geq
            \pi\sqrt{p(p+1)}
            - \pi\sqrt{p(p+1)} p(p+1) \sum_{n=-\infty}^\infty |g_n|^2 -C_p\delta^{1/2}\Vert g\Vert_{H^{1/2}}^2 \\
            &\quad + \frac{p}{2p+1} \pi\sum_{n=-\infty}^{\infty} (-\operatorname{sgn}(n)) \sqrt{n(n+1)} \vert F^-_n \vert^2
            + \frac{p+1}{2p+1}\pi \sum_{n=-\infty}^{\infty}\sqrt{n(n+1)}\vert F^+_n\vert^2.
        \end{aligned}
    \end{equation}
    In the rest of the proof, we will use the high-frequency terms in the second line to overcome the negative low-frequency terms in the sum of the first line.
    
    Let us introduce the quadratic errors
    $$R^+(e^{i\eta}):=F^+(e^{i\eta})-ipe^{ip\eta}g(\eta),
    \qquad R^-(e^{i\eta}):=F^-(e^{i\eta})+i(p+1)e^{-i(p+1)\eta}g(\eta).$$
    Clearly their Fourier coefficients satisfy
    \[
        F^+_n = R^+_n + ip g_{n-p}
        \qquad \text{and} \qquad
        F^-_n = R^-_n - i(p+1) g_{n+p+1}.
    \]
    We claim that
    \begin{align}
        \Vert R^+\Vert_{H^{1/2}}+\Vert R^-\Vert_{H^{1/2}}\leq C_p\Vert g\Vert_{L^{\infty}}\Vert g\Vert_{H^{1/2}}.
        \label{eq:estimates_R+_R-}
    \end{align}
    Indeed, note that $R^+(e^{i\eta}) = e^{ip\eta}r(g(\eta))$, where $r(t) := e^{ipt}-1-ipt$. The function $r$ satisfies
    \[
    \lvert r(t)\rvert\leq \frac{p^2}{2}t^2,\quad\lvert r'(t)\rvert\leq p^2\lVert g\rVert_{L^\infty} \quad\text{ for }t\in[-\lVert g\rVert_{L^\infty}, \lVert g\rVert_{L^\infty}].
    \]
    Thus, for any $\eta, \kappa \in \mathbb R$,
    \[
    \lvert r(g(\eta))\rvert\leq \frac{p^2}{2}\lvert g(\eta)\rvert^2,
    \qquad \lvert r(g(\eta))-r(g(\kappa))\rvert\leq p^2\lVert g\rVert_\infty\lvert g(\eta)-g(\kappa)\rvert. 
    \]
    Working with the difference quotient seminorm of $H^{1/2}$, one deduces that $r\circ g\in H^{1/2}$ with $\lVert r\circ g\rVert_{H^{1/2}}\leq C_p\lVert g\rVert_\infty\lVert g\rVert_{H^{1/2}}$. Since the multiplication by $e^{ip\eta}$ is continuous in $H^{1/2}$, we deduce that $R^+$ satisfies \eqref{eq:estimates_R+_R-}. A similar argument shows that $R^-$ satisfies \eqref{eq:estimates_R+_R-} as well.

    Therefore, as $\sqrt{n(n+1)}\leq C(1+\vert n\vert)$, by expanding the squares and using \eqref{eq:estimates_R+_R-} and Cauchy--Schwarz in $H^{1/2}$, we get
    \begin{equation}
        \sum_{n \in \Z} \sqrt{n(n+1)}\vert F^+_n\vert^2\geq p^2\sum_{n \in \Z}\sqrt{(n+p)(n+p+1)}\vert g_n\vert^2-C_p \Vert g\Vert_{L^{\infty}} \left( 1 + \Vert g\Vert_{L^{\infty}}  \right) \Vert g\Vert_{H^{1/2}}^2\label{eq: turning F to g, positive}
    \end{equation}
    and 
    \begin{equation}
        \label{eq: turning F to g, negative}
        \begin{aligned}
            -\sum_{n\in\mathbb Z}\operatorname{sgn}(n)
            \sqrt{n(n+1)}|F^-_n|^2\geq&
            -(p+1)^2\sum_{n\in\mathbb Z}
            \operatorname{sgn}(n-p-1)
            \sqrt{(n-p-1)(n-p)}\,|g_n|^2\\
            &- C_p \|g\|_{L^\infty} (1 + \|g\|_{L^\infty}) \|g\|_{H^{1/2}}^2.
        \end{aligned}
    \end{equation}
    
    Recall that by our choice of $g$, $g_n=0$ for $\vert n\vert\leq p$. Also, since $g$ is real-valued, we have $g_n=\overline{g_{-n}}$.
    Substituting \eqref{eq: turning F to g, positive} and \eqref{eq: turning F to g, negative} into \eqref{eq: aux lower bound E} and reindexing to only sum over positive values, we get
    \begin{equation}
        \label{eq:estimate-EG-minus-cone}
        E(u) - \pi\sqrt{p(p+1)}
        \geq \pi p(p+1) \sum_{n=p+1}^\infty A_n |g_n|^2 - C_p\delta^{1/2}\Vert g\Vert_{H^{1/2}}^2,
    \end{equation}
    where
    \[
        A_n := \sqrt{(n+p)(n+p+1)} - \frac{1}{2p+1} \sqrt{(n-p)(n-p-1)} - 2 \sqrt{p(p+1)}.
    \]
    The crucial inequality\footnote{Instead, $A_n < 0$ for $0 \leq n \leq p$.} is
    \begin{equation}
        \label{eq: critical n}
        A_n \geq \frac{2p}{2p+1} \left( n - p - \frac{1}{2} \right) \geq c_p (n + 1)
        \qquad \text{ for all } n \geq p+1,
    \end{equation}
    where $c_p > 0$ is some constant. From this, \eqref{eq:estimate-EG-minus-cone} implies that
    $$E(u) - \pi\sqrt{p(p+1)} \geq (c_p - C_p\delta^{1/2}) \Vert g\Vert_{H^{1/2}}^2 \geq 0,$$
    provided that $\delta>0$ is chosen sufficiently small.

    Finally, to prove \eqref{eq: critical n}, define $f(t) := \sqrt{t(t+1)}$ and observe that $f(t) < t + \frac{1}{2}$ for $t \geq 0$ and $f'(t) > 1$. By the first property, $f(n-p-1) + f(p) < n$, and by the second one, $f(n+p) > f(p) + n$. Hence
    \[
        f(n+p) - \frac{1}{2p+1} f(n-p-1) - 2 f(p)
        > f(p) + n - \frac{1}{2p+1} (n - f(p)) - 2 f(p)
        = \frac{2p (n - f(p))}{2p+1}
    \]
    and the conclusion follows from $f(p) < p + \frac{1}{2}$.
\end{proof}

\begin{lem}\label{lem: lower bound for perturbed single pole problem}

    Let $\varphi:\Sp^1\rightarrow\Sp^1$ be a $C^1$ diffeomorphism satisfying $\lVert \varphi-\Id_{\Sp^1}\rVert_{C^1}<\delta$ and let $z_0\in \Disk$ with $\vert z_0\vert<\delta$, for $\delta > 0$ small. Then for any $F\in C^1(\Sp^1,\CC)$ and $v\in H^1(\Disk,\CC)$ with $v|_{\p\Disk}=F\circ\varphi$, there holds
    $$
    E_{G_{z_0}}(v)
    = E(v)+\int_{\Disk}G_{z_0}J(v)\dif x\dif y
    \geq \pi\sum_{n=-\infty}^{\infty}\sqrt{n(n+1)}\vert F_n\vert^2-C\delta^{1/2}\Vert F\Vert_{H^{1/2}}^2,$$
    where $C>0$ and the Fourier coefficients and norms of $F$ are taken with respect to the coordinate $\eta$ induced by $\varphi$.
\end{lem}
\begin{proof}
    Let $b(z)=(z-z_0)/(1-\bar z_0z)$ be the biholomorphic diffeomorphism of $\Disk$ which maps $z_0$ to $0$, and let $\zeta$ be the coordinate on $\Sp^1$ induced by $b$, so that $b(e^{i\theta})=e^{i\zeta}$. By Lemma \ref{lem: minimizing Fourier expansion of energy}, there holds
    $$ E_{G_{z_0}}(v) = E_{G_{0}}(v \circ b^{-1}) \geq \pi\sum_{n=-\infty}^{\infty}\sqrt{n(n+1)}\vert \tilde F_n \vert^2,$$
    where $\tilde F := F \circ \varphi \circ b^{-1} : \Sp^1_\zeta \to \CC$ and its Fourier coefficients are taken with respect to $\zeta$.  
    The Gagliardo seminorm of $F : \Sp^1_\eta \to \CC$, with respect to the coordinate $\eta$, is
    $$2\pi \sum_{n=-\infty}^{\infty}\vert n\vert \vert F_n\vert^2=[F]_{\dot H^{1/2}(\dif\eta)}^2=\frac{1}{2\pi}\int_0^{2\pi}\int_{0}^{2\pi}\frac{\vert F(e^{i\eta})-F(e^{i\kappa})\vert^2}{\vert e^{i\eta}-e^{i\kappa}\vert^2}\,\dif \eta\,\dif \kappa.$$
    By assumption, we can write $e^{i\eta(\zeta)} = \varphi\circ b^{-1}(e^{i\zeta})=e^{i(\zeta+f(\zeta))}$ with $\Vert f\Vert_{C^1}\leq C\delta$.
    If $h \in C^1(\mathbb R, \mathbb R)$ denotes the inverse of $\zeta \mapsto \zeta + f(\zeta)$, so that $\zeta = h(\eta)$, the analogous expression for $\tilde F$ and the change of variables formula give
    \begin{align*}
        2\pi \sum_{n=-\infty}^{\infty}\vert n\vert\vert \tilde F_n\vert^2
        &=\frac{1}{2\pi}\int_0^{2\pi}\int_{0}^{2\pi}\frac{\vert \tilde F(e^{i\zeta})-\tilde F(e^{i\xi})\vert^2}{\vert e^{i\zeta}-e^{i\xi}\vert^2}\,\dif \zeta\,\dif \xi \\
        &=\frac{1}{2\pi}\int_0^{2\pi}\int_{0}^{2\pi}\frac{\vert F(e^{i\eta})-F(e^{i\kappa})\vert^2}{\vert e^{i\eta}-e^{i\kappa}\vert^2} K(\eta, \kappa) \,\dif \eta\,\dif \kappa,
    \end{align*}
    where
    \[
        K(\eta, \kappa) := \frac{\left| e^{i\eta} - e^{i\kappa} \right|^2}{\left| e^{ih(\eta)} - e^{ih(\kappa)} \right|^2} h'(\eta) h'(\kappa)
        \quad \text{satisfies} \quad
        \left\| K - 1 \right\|_{L^\infty} \leq C \delta.
    \]
    This follows from $\vert h'-1\vert\leq C\delta$ and $\vert h(\eta)-\eta\vert\leq C\delta$. Hence,
    $$\left\vert\sum_{n=-\infty}^{\infty}\vert n\vert(\vert \tilde F_n\vert^2-\vert F_n\vert^2)\right\vert\leq C\delta \| F \|_{H^{1/2}}^2.$$
    Since $\left| \sqrt{n(n+1)} - |n|\right| \leq 1 \, \forall n$, to complete the proof it will suffice to show the estimate
    \[
        \sum_{n=-\infty}^{\infty} \left\lvert \vert \tilde F_n\vert^2-\vert F_n\vert^2 \right\rvert
        \leq \sum_{n=-\infty}^{\infty}
        \lvert \tilde F_n - F_n \rvert
        ( \vert \tilde F_n\vert+\vert F_n\vert )
        \leq C\Vert \tilde F-F\Vert_{L^2} (\Vert F\Vert_{L^2}+\Vert \tilde F\Vert_{L^2} )\leq C\delta^{1/2}\Vert F\Vert_{H^{1/2}}^2.
    \]
    The bound on the second factor is clear, again by change of variables formula:
    \begin{equation}
        \label{eq: L2 F tilde F}
        \lVert \tilde F\rVert_{L^2}^2=\int_0^{2\pi}\lvert F(e^{i\eta})\rvert^2h'(\eta) \,\dif \eta\leq C\lVert F\rVert_{L^2}^2.
    \end{equation}
    It only remains to show that the operator
    $T F(e^{i\zeta}) := F(e^{i(\zeta + f(\zeta))}) - F(e^{i\zeta})$
    satisfies the bound
    $\| T F\|_{L^2} \leq C \delta^{1/2} \| F \|_{H^{1/2}}$ for all $F \in H^{1/2}$. By interpolation, this follows from the estimates $\| T \|_{L^2 \to L^2} \leq C$ and $\| T \|_{H^{1} \to L^2} \leq C \delta$; the first one is immediate from~\eqref{eq: L2 F tilde F} and the triangle inequality, and the second one follows from the following computation: write
    \[
        T F(e^{i\zeta})
        = \tilde F(e^{i\zeta})-F(e^{i\zeta})
        = f(\zeta)\int_0^1(\partial_\theta F)(e^{i(\zeta+tf(\zeta))}) \, \dif t, \quad \text{with }\lVert f\rVert_{C^1}\leq C\delta,
    \]
    and denote by $h_t$ the inverse of $\zeta \mapsto \zeta + t f(\zeta)$, so that
    \[
        \int_{0}^{2\pi} |T F(e^{i\zeta})|^2 \, \dif \zeta
        \leq C \delta^2 \int_{0}^{2\pi} \int_0^1 \left| (\partial_\theta F) (e^{i(\zeta+tf(\zeta))}) \right|^2 \, \dif t \, \dif \zeta
        \leq C \delta^2 \int_0^1 \int_{0}^{2\pi} \left| (\partial_\theta F) (e^{i\kappa}) \right|^2 \, h_t'(\kappa) \dif \kappa \, \dif t
    \]
    and the claim follows since $|h_t'| \leq C$.
\end{proof}

\bibliography{biblio}
\bibliographystyle{plain}
\end{document}